\documentclass[11pt]{amsart}
\usepackage{amssymb, adjustbox, enumerate, amsbsy, stmaryrd}
\usepackage{geometry}

\usepackage{todonotes}

\usepackage{amsfonts, amssymb, amscd}
\numberwithin{equation}{section}

\usepackage[symbol]{footmisc}

\usepackage{bm}
\usepackage{verbatim}
\usepackage{mathrsfs}
\usepackage{graphicx}
\usepackage{tikz-cd}
\usepackage{subcaption}
\usepackage{listings}
\usepackage{subfiles}
\usepackage[toc,page]{appendix}
\usepackage{mathtools}
\usepackage{comment}
\usepackage{enumerate}
\usepackage{enumitem}
\usepackage[all]{xy}

\usepackage{graphicx}
\graphicspath{{images/}}

\usepackage{appendix}
\usepackage{hyperref}
\hypersetup{
    colorlinks=true,
    citecolor=red,
    linkcolor=blue,
    filecolor=magenta,      
    urlcolor=red,
}
\newcommand{\Qq}{\mathbb{Q}}

\newcommand{\Zz}{\mathbb{Z}}

\newcommand{\mld}{{\rm{mld}}}

\newcommand{\Ii}{\Gamma}

\newcommand{\mm}{\mathfrak{m}}

\newcommand\coeff{{\rm{coeff}}}

\newtheorem{thm}{Theorem}[section]

\newtheorem{lem}[thm]{Lemma}
\newtheorem{sett}[thm]{Setting}
\newtheorem{prop}[thm]{Proposition}

\newtheorem{claim}[thm]{Claim}

\theoremstyle{definition}
\newtheorem{defn}[thm]{Definition}

\theoremstyle{definition}
\newtheorem{rem}[thm]{Remark}

\theoremstyle{definition}

\begin{document}

\title{Classification of threefold enc cDV quotient singularities}
\author{Jingjun Han and Jihao Liu}

\subjclass[2020]{14E30,14C20,14E05,14J17,14J30,14J35}
\date{\today}

\address{Shanghai Center for Mathematical Sciences, Fudan University, Jiangwan Campus, Shanghai, 200438, China}
\email{hanjingjun@fudan.edu.cn}

\address{Department of Mathematics, Peking University, No. 5 Yiheyuan Road, Haidian District, Peking 100871, China}
\email{liujihao@math.pku.edu.cn}

\begin{abstract}
    We provide a rough classification of threefold exceptionally non-canonical cDV quotient singularities by studying their combinatorial behavior.
\end{abstract}

\maketitle

\tableofcontents

\section{Introduction}

We work over the field of complex numbers $\mathbb C$.

This is the note to the paper ``On termination of flips and exceptionally non-canonical singularities" \cite{HL22} by the authors. This note appeared as the appendix of the arXiv version of the paper \url{https://arxiv.org/pdf/2209.13122.pdf} but is not going to be included in the final published version of \cite{HL22}. This is because the proof of the main theorem of this note is elementary but requires complicated combinatorical computations. By suggestions from the referee(s) and for the reader's convenience, we take the appendix of \cite{HL22} out and write it as this separate, self-contained note.

The goal of this note is to provide a (rough) classification of enc cDV (cyclic) quotient singularities. We start with the following setting.

\begin{sett}\label{Setting: before terminal lem}
We set up the following notations and conditions.
\begin{enumerate}
\item Let $r$ be a positive integer, $0\leq a_1,a_2,a_3,a_4,e<r$ integers, such that \begin{enumerate}
       \item $\gcd(a_i,r)\mid\gcd(e,r)$ for any $1\leq i\leq 4$.
    \item $\gcd(a_i,a_j,r)=1$ for any $1\le i<j\le 4$.
    \item  $\sum_{i=1}^4a_i-e\equiv 1\mod r$.
\end{enumerate}
\item $f\in \mathbb C\{x_1,x_2,x_3,x_4\}$ is $\bm{\mu}$-semi-invariant, that is, ${\bm{\mu}}(f)=\xi^ef$, and is one of the following $3$ types:
\begin{enumerate}
    \item (cA type) $f=x_1x_2+g(x_3,x_4)$ with $g\in\mm^2$.
    \item (Odd type) $f=x_1^2+x_2^2+g(x_3,x_4)$ with $g\in\mm^3$ and $a_1\not\equiv a_2\mod r$.
    \item (cD-E type) $f=x_1^2+g(x_2,x_3,x_4)$ with $g\in\mm^3$,
\end{enumerate}
where $\mm$ is the maximal ideal of $\mathbb C\{x_1,x_2,x_3,x_4\}$, and ${\bm{\mu}}:\mathbb C^4\rightarrow\mathbb C^4$ is the action $(x_1,x_2,x_3,x_4)\rightarrow (\xi^{a_1}x_1,\xi^{a_2}x_2,\xi^{a_3}x_3,\xi^{a_4}x_4)$.
\item 
One of the two cases hold:
\begin{enumerate}
    \item $\alpha(x_1x_2x_3x_4)-\alpha(f)>1$ for any $\alpha\in N$. In this case, we let $k:=1$ and $\beta:=\bm{0}$.
    \item There exists an integer $k\geq 2$, and a primitive vector $\beta\in N$, such that either
    \begin{enumerate}
        \item 
        \begin{itemize}
            \item $\frac{1}{k}<\beta(x_1x_2x_3x_4)-\beta(f)\leq\min\{\frac{12}{13},\frac{1}{k-1}\}$, or
            \item $\beta(x_1x_2x_3x_4)-\beta(f)=1$ and $k=2$, 
        \end{itemize} 
        and
        \item for any $\alpha\in N\backslash\{\beta,2\beta,\dots,(k-1)\beta\}$, $\alpha(x_1x_2x_3x_4)-\alpha(f)>1$,
    \end{enumerate}
\end{enumerate}
where $$N:=\{w\in\mathbb Q^4_{\geq 0}\mid w\equiv\frac{1}{r}(ja_1,ja_2,ja_3,ja_4)\mod \mathbb Z^4\text{ for some }j\in\mathbb Z\}\backslash\{\bm{0}\}.$$
\end{enumerate}
Moreover, if $f$ is of cA type, then for any integer $a$ such that $\gcd(a,r)=1$, $\frac{1}{r}(a_1,a_2,a_3,a_4,e)\not\equiv\frac{1}{r}(a,-a,1,0,0)\mod \Zz^5$. 
\end{sett}

The main theorem of this note is the following:

\begin{thm}\label{thm: note beta finite}
Notations and conditions as in Setting \ref{Setting: before terminal lem}. Then either $r$ or $\beta\not=\bm{0}$ belongs to a finite set depending only on $k$.
\end{thm}

Theorem \ref{thm: note beta finite} implies the following theorem and we refer to \cite[Theorem 6.8]{HL22} for a proof:

\begin{thm}\label{thm: enc cDV cover case}
Let $\Ii\subset[0,1]$ be a DCC (resp. finite) set. Assume that $(X\ni x,B)$ is a $\Qq$-factorial enc pair of dimension $3$, such that
\begin{enumerate}
    \item $X\ni x$ is an isolated non-canonical singularity, 
    \item $\coeff(B)\subseteq \Ii$, and
    \item $\tilde X\ni\tilde x$ is terminal but not smooth, where $\pi: (\tilde X\ni\tilde x)\rightarrow (X\ni x)$ is the index $1$ cover of $X\ni x$.
\end{enumerate}
Then $\mld(X\ni x,B)$ belongs to an ACC set (resp. is discrete away from $0$).
\end{thm}

\begin{rem}
If we could show that $r$ belongs to a finite set, then $a_i,e$ also belongs to a finite set, hence the action of ${\bm{\mu}}$ on $\mathbb C^4$ and $f$ belongs to a finite set. Since the singularities $X\ni x$ in Theorem \ref{thm: enc cDV cover case} are the cDV quotient singularities of type
 $$(\mathbb C^4\supset (f=0))/{\bm{\mu}},$$
it could be regarded as that we classify all such kinds of singularities. This is why we say that Theorem \ref{thm: note beta finite} gives a (rough) classification of enc cDV quotient singularities. One difficulty at the moment is that we can only show that ``either $r$ or $\beta\not=\bm{0}$ belongs to a finite set" and could not show that $r$ belongs to a finite set.
\end{rem}

Theorem \ref{thm: note beta finite} is a consequence of Theorems \ref{thm: cA type beta finite} and \ref{thm: non cA type beta finite} below which will be proven in Sections \ref{note: cA type} and \ref{note: noncA type} respectively.

\begin{thm}\label{thm: cA type beta finite}
Notations and conditions as in Setting \ref{Setting: before terminal lem}. For each positive integer $k$, there exists a finite set $\Ii_k'$ depending only on $k$ satisfying the following. If $f$ is of cA type, then either $r\in \Ii'_k$ or $\bm{0}\neq\beta\in \Ii'_k$.
\end{thm}

\begin{thm}\label{thm: non cA type beta finite}
Notations and conditions as in Setting \ref{Setting: before terminal lem}. For each positive integer $k$, there exists a finite set $\Ii_k'$ depending only on $k$ satisfying the following. If $f$ is not of cA type, then either $r\in \Ii'_k$ or $\bm{0}\neq\beta\in \Ii'_k$.
\end{thm}

\noindent\textbf{Acknowledgement.} The authors would like to thank the referee(s) of \cite{HL22} for useful suggestions. The work is supported by the National Key R\&D Program of China (\#2024YFA1014400, \#2023YFA1010600, \#2020YFA0713200). The first author is supported by NSFC for Excellent Young Scientists (\#12322102).

\section{Preliminaries}

We will use the following definition in the rest part of this note.

\begin{defn}
Settings as in Setting \ref{Setting: before terminal lem}. Let
\begin{itemize}
\item $N^0:=N\cap [0,1]^4\backslash\{0,1\}^4$.
    \item $\alpha_j:=(\{\frac{ja_1}{r}\},\{\frac{ja_2}{r}\},\{\frac{ja_3}{r}\},\{\frac{ja_4}{r}\})$ for any $1\leq j\leq r-1$, and
    \item $w':=(1,1,1,1)-w$ for any $w\in\mathbb Q^d$,
\end{itemize} 
We define three sets $\Psi_1,\Psi_2$ and $\Psi$ in the following way. If $k=1$, then let $\Psi_1:=\Psi_2:=\Psi:=\emptyset$. If $k\geq 2$, then we let
\begin{itemize}
    \item $\Psi_1:=\{\beta,2\beta,\dots,(k-1)\beta\}$, 
    \item $\Psi_2:=\{\beta',(2\beta)',\dots,((k-1)\beta)'\}$, and
    \item $\Psi:=\Psi_1\cup\Psi_2$.
\end{itemize}
\end{defn}

The following lemma appeared in \cite{HL22}. For the reader's convenience and in order to make this note self-contained, we provide its full proof here.

\begin{lem}[{\cite[Lemma 6.2]{HL22}}]\label{lem: toric bdd index lemma}
Let $d$ be a positive integer and $\epsilon$ a positive real number. Then there exists a positive integer $I$, depending only on $d$ and $\epsilon$, satisfying the following. Let $r$ be a positive integer and $v_1,\dots,v_d\in [0,1]$ real numbers, such that  $\sum_{i=1}^d(1+(m-1)v_i-\lceil mv_i\rceil)\ge\epsilon$ for any $m\in [2,r]\cap \Zz$. Then $r\le I$.
\end{lem}
\begin{proof}
Suppose that the statement does not hold. Then for each $j\in\Zz_{\ge 1}$, there exist $v_{1,j},\dots,v_{d,j}\in [0,1]$ and positive integers $r_j$, such that
\begin{itemize}
    \item $\sum_{i=1}^d(1+(m-1)v_{i,j}-\lceil mv_{i,j}\rceil)\ge\epsilon$
for any $m\in [2,r_j]\cap\Zz$,
\item $r_j$ is strictly increasing, and
\item $\bar v_i:=\lim_{j\rightarrow+\infty}v_{i,j}$ exists.
\end{itemize}
Let $\bm{v}:=(\bar v_1,\dots,\bar v_d)$. 
By Kronecker's theorem, there exist a positive integer $n$ and a vector $\bm{u}\in\mathbb Z^{d}$ such that $||n\bm{v}-\bm{u}||_{\infty}<\min\{\frac{\epsilon}{d},\bar v_i\mid \bar v_i>0\}$ and $n\bar v_i\in\Zz$ for any $i$ such that $\bar v_i\in\mathbb Q$. In particular, $\lceil (n+1)\bar v_i\rceil=\lfloor (n+1)\bar v_i\rfloor+1$ for any $i$ such that $\bar v_i\in (0,1)$. Now $\lim_{j\rightarrow+\infty}(1+nv_{i,j}-\lceil (n+1)v_{i,j}\rceil)=0$ when $\bar v_i=0$ and $\lim_{j\rightarrow+\infty}(1+nv_{i,j}-\lceil (n+1)v_{i,j}\rceil)=1+n\bar v_i-\lceil (n+1)\bar v_i\rceil$ when $\bar v_i>0$. Thus
\begin{align*}
    &\lim_{j\rightarrow+\infty}\sum_{i=1}^d(1+nv_{i,j}-\lceil (n+1)v_{i,j}\rceil)=\sum_{0<\bar v_i<1}(1+n\bar v_i-\lceil (n+1)\bar v_i\rceil)\\
    =&\sum_{0<\bar v_i<1}(1+(n+1)\bar{v_i}-\lceil (n+1)\bar v_i\rceil-\bar v_i)=\sum_{\bar 0<\bar v_i<1}(\{(n+1)\bar v_i\}-\bar v_i)<\sum_{0<\bar v_i<1}\frac{\epsilon}{d}\le \epsilon.
\end{align*}
Thus possibly passing to a subsequence, $\sum_{i=1}^d(1+nv_{i,j}-\lceil (n+1)v_{i,j}\rceil)<\epsilon$ for any $j$, hence $n>r_j$, which contradicts $\lim_{j\rightarrow+\infty}r_j=+\infty$.
\end{proof}

Theorem \ref{thm: terminal lemma} is known as the terminal lemma.
\begin{thm}[{cf. \cite[(5.4) Theorem, (5.6) Corollary]{Rei87}, \cite[Theorem 2.6]{Jia21}}] \label{thm: terminal lemma}
Let $r$ be a positive integer and $a_1,a_2,a_3,a_4,e$ integers, such that $\gcd(a_4,r)=\gcd(e,r)$ and $\gcd(a_1,r)=\gcd(a_2,r)=\gcd(a_3,r)=1$. Suppose that
$$\sum_{i=1}^4\{\frac{ja_i}{r}\}=\{\frac{je}{r}\}+\frac{j}{r}+1$$
for any integer $1\leq j\leq r-1$. We have the following.
\begin{enumerate}
    \item If $\gcd(e,r)>1$, then $a_4\equiv e\mod r$, and there exists $i_1,i_2,i_3$ such that $\{i_1,i_2,i_3\}=\{1,2,3\}$, $a_{i_1}\equiv 1\mod r$, and $a_{i_2}+a_{i_3}\equiv 0\mod r$.
    \item If $\gcd(e,r)=1$, we let $a_5:=-e$ and $a_6:=-1$, then there exists $i_1,i_2,i_3,i_4,i_5,i_6$ such that $\{i_1,i_2,i_3,i_4,i_5,i_6\}=\{1,2,3,4,5,6\}$, such that $a_{i_1}+a_{i_2}\equiv a_{i_3}+a_{i_4}\equiv a_{i_5}+a_{i_6}\equiv 0\mod r$.
\end{enumerate}
\end{thm}

In order to study enc singularities, Jiang introduced the so-called non-canonical lemma \cite[Lemma 2.7]{Jia21}. The following lemma is a generalization of his result. 

\begin{lem}\label{lem: transfer to fivefold lemma refined}
Let $\delta$ be a positive real number. Then there exists a finite set $\Ii_0\subset (0,1)\cap\mathbb Q$ depending only on $\delta$ satisfying the following. Assume that $k_0,r$ are two positive integers such that $1\leq k_0\leq r-1$, and $a_1,a_2,a_3,a_4,e$ five integers, such that
\begin{enumerate}
    \item $$\sum_{i=1}^4\{\frac{a_ik_0}{r}\}=\{\frac{ek_0}{r}\}+\frac{k_0}{r},$$
    and
    \item for any $k\not=k_0$ such that $1\leq k\leq r-1$, $$\sum_{i=1}^4\{\frac{a_ik}{r}\}\geq \{\frac{ek}{r}\}+\frac{k_0}{r}+\delta.$$
\end{enumerate}
Then $\frac{r}{\gcd(r,k_0)}\in\Ii_0$. In particular, $\frac{k_0}{r}$ also belongs to a finite set depending only on $\delta$.
\end{lem}

\begin{proof}
Let $v_i:=\{\frac{k_0a_i}{r}\}$ for any $1\leq i\leq 4$, and $v_5:=\{\frac{(r-e)k_0}{r}\}$. For any $2\le m\le \frac{r}{\gcd(r,k_0)}-1$, $mk_0\equiv k \mod r$ for some $1\le k\le r-1$, then $k\neq k_0$ as $r\nmid (m-1)k_0$. For any $0\le a\le 1$, we have
$$1+ m\{1-a\}-\lceil m\{1-a\} \rceil=1-\{ma\}. $$
In particular,
$$1+mv_5-\lceil mv_5 \rceil=1-\{\frac{emk_0}{r}\}.$$
Moreover, for any $v\ge 0$, we have $1+v-\lceil v\rceil \ge \{v\}$. Thus
\begin{align*}\sum_{i=1}^ 5(1+(m-1)v_i-\lceil mv_i\rceil)
=&\sum_{i=1}^5(1+mv_i-\lceil mv_i\rceil)-\sum_{i=1}^5 v_i\\
\ge &\sum_{i=1}^4\{mv_i\}-\sum_{i=1}^4 v_i+(1+mv_5-\lceil mv_5\rceil)-(1+v_5-\lceil v_5\rceil)\\
=&\sum_{i=1}^4\{mv_i\}-\{\frac{emk_0}{r}\}-(\sum_{i=1}^4 \{v_i\}-\{\frac{ek_0}{r}\})  \\
\ge &1+\frac{k_0}{r}+\delta-1-\frac{k_0}{r}=\delta.
\end{align*}
By Lemma \ref{lem: toric bdd index lemma}, $\frac{r}{\gcd(r,k_0)}\in \Ii_0$ for some set $\Ii_0$ depending only on $\delta$.
\end{proof}

\section{cA type}\label{note: cA type}

The goal of this subsection is to show Theorem \ref{thm: cA type beta finite}. 


\begin{thm}\label{thm: cA case up to terminal lemma}
Notations and conditions as in Setting \ref{Setting: before terminal lem}. For each positive integer $k$, there exists a finite set $\Ii_k\subset\mathbb Q^4\cap [0,1]^4$ depending only on $k$ satisfying the following. Suppose that $f$ is of cA type. 
Then possibly switching $x_3$ and $x_4$, either $\bm{0}\not=\beta\in\Ii_k$, or the following holds.
\begin{enumerate}
    \item For any $\alpha\in N^0\backslash\Psi$, there exists $w\in\{\alpha,\alpha'\}$, such that
    \begin{enumerate}
        \item $w(f)=w(x_1x_2)\leq 1$ and $w'(f)=w(x_1x_2)-1$, 
        \item $w(x_3x_4)>1$ and $w'(x_3x_4)<1$, and
        \item $w(x_1x_2)=1$ if and only if $w'(x_1x_2)=1$. Moreover, if  $w(x_1x_2)=1$, then either $w(x_3)=1$ or $w(x_4)=1$, and either $w'(x_3)=0$ or $w'(x_4)=0$.
    \end{enumerate}
    \item For any $1\leq j\leq r-1$ such that $\alpha_j\not\in\Psi$, either
    \begin{itemize}
        \item  $\{\frac{ja_1}{r}\}+\{\frac{ja_2}{r}\}=\{\frac{je}{r}\}$ and $\{\frac{ja_3}{r}\}+\{\frac{ja_4}{r}\}=\frac{j}{r}+1$, or
        \item $\{\frac{ja_1}{r}\}+\{\frac{ja_2}{r}\}=\{\frac{je}{r}\}+1$ and $\{\frac{ja_3}{r}\}+\{\frac{ja_4}{r}\}=\frac{j}{r}$.
    \end{itemize}
    \item $\gcd(a_1,r)=\gcd(a_2,r)=\gcd(a_3,r)=1$ and $\gcd(a_4,r)=\gcd(e,r)$.
    \item If $\beta\in N^0$, then there exists $1\leq k_0\leq r-1$, such that 
    \begin{enumerate}
        \item $\beta=\alpha_{k_0}$,
        \item $\frac{1}{k}<\beta(x_3x_4)=\frac{k_0}{r}<\min\{\frac{13}{14},\frac{1}{k-1}\}$, and
            \item $\beta(x_1x_2)\geq 1$, $\{\frac{k_0a_1}{r}\}+\{\frac{k_0a_2}{r}\}=\{\frac{k_0e}{r}\}+1$, and $\{\frac{k_0a_3}{r}\}+\{\frac{k_0a_4}{r}\}=\frac{k_0}{r}$.
    \end{enumerate}
    \item  For any $1\leq j\leq r-1$,
    $$\sum_{i=1}^4\{\frac{ja_i}{r}\}=\{\frac{je}{r}\}+\frac{j}{r}+1.$$
\end{enumerate}
\end{thm}
\begin{proof}
\noindent\textbf{Step 1}. In this step we summarize some auxiliary results that will be used later.

Since $x_1x_2\in f$ and $\sum_{i=1}^4 a_i-e\equiv 1 \mod r$, we have $a_1+a_2\equiv e\mod r$, and $a_3+a_4\equiv 1\mod r$. By Setting \ref{Setting: before terminal lem}(1.a)(1.b), $\gcd(a_1,r)=\gcd(a_2,r)=1$. 

Since $a_1+a_2\equiv e\mod r$ and ${\bm{\mu}}(f)=\xi^ef$, $\alpha(f)\equiv \alpha(x_1x_2)\mod \Zz$ for any $\alpha\in N$. 

For any $\alpha\in N^0$, by Setting \ref{Setting: before terminal lem}(1.b), $0<\alpha(x_1x_2)<2$. Since $x_1x_2\in f$, $0\leq\alpha(f)\leq\alpha(x_1x_2)<2$, hence either $\alpha(f)=\alpha(x_1x_2)$ or $\alpha(f)=\alpha(x_1x_2)-1$. By Setting \ref{Setting: before terminal lem}(3.b), for any $\alpha\in N^0\backslash\Psi_1$ such that $\alpha(f)=\alpha(x_1x_2),\alpha(x_3x_4)>1$.

We may assume that $k$ is fixed. If $\beta\neq \bm{0}$, then since $\frac{1}{k}<\beta(x_1x_2x_3x_4)-\beta(f)\leq\frac{1}{k-1}$, $\beta(x_1x_2)\equiv\beta(f)\mod \Zz$, and $\beta(f)\leq\beta (x_1x_2)$, we have $\beta(x_1x_2)=\beta(f)$ and $\frac{1}{k}<\beta(x_3x_4)\leq\frac{1}{k-1}$. 

Finally, since switching $x_3$ and $x_4$ will not influence (1)(2)(4)(5), we will only have a possibly switching of $x_3$ and $x_4$ only when we prove (3).

\medskip

\noindent\textbf{Step 2}. In this step we prove (1). Pick $\alpha\in N^0\backslash\Psi$, then $\alpha'\in N^0\backslash\Psi$. Since $0<\alpha(x_1x_2)<2$ and $0<\alpha'(x_1x_2)<2$, $\alpha(x_1x_2)=1$ if and only if $\alpha'(x_1x_2)=1$. By \textbf{Step 1}, there are two cases:

\medskip

\noindent\textbf{Case 1}. $\alpha(f)=\alpha(x_1x_2)-1$. In this case, $\alpha(x_1x_2)\geq 1$. There are two sub-cases:

\medskip

\noindent\textbf{Case 1.1}. If $\alpha(x_1x_2)=1$, then $\alpha(f)=0$. Since $\gcd(a_1,r)=\gcd(a_2,r)=1$, $\alpha(x_1)\not=0$, and $\alpha(x_2)\not=0$. Thus either $\alpha(x_3)=0$ or $\alpha(x_4)=0$, and either $\alpha'(x_3)=1$ or $\alpha'(x_4)=1$. By Setting \ref{Setting: before terminal lem}(1.b), $\alpha(x_3x_4)<1$.

\medskip

\noindent\textbf{Case 1.2}. If $\alpha(x_1x_2)>1$, then $\alpha'(x_1x_2)<1$, hence $\alpha'(f)=\alpha'(x_1x_2)$. By Setting \ref{Setting: before terminal lem}(3.b), $\alpha'(x_3x_4)>1$, hence $\alpha(x_3x_4)<1$.

\medskip

In either sub-case, $\alpha(x_3x_4)<1$, hence $\alpha'(x_3x_4)>1$. Therefore, we may take $w=\alpha'$. Moreover, $\alpha'$ is not of \textbf{Case 1} as $\alpha'(x_1x_2)\le 1$, hence $\alpha'(f)=\alpha'(x_1x_2)$.

\medskip

\noindent\textbf{Case 2}. $\alpha(f)=\alpha(x_1x_2)$. In this case, by \textbf{Step 1}, $\alpha(x_3x_4)>1$, so $\alpha'(x_3x_4)<1$, $\alpha'(f)\not=\alpha'(x_1x_2)$, and $\alpha'(f)=\alpha'(x_1x_2)-1$. Thus $\alpha'$ is of \textbf{Case 1}. By \textbf{Case 1}, $\alpha(x_1x_2)\leq 1$. 

Moreover, if $\alpha(x_1x_2)=1$, then since $\alpha'$ is of \textbf{Case 1}, either $\alpha'(x_3)=0$ or $\alpha'(x_4)=0$, hence either $\alpha(x_3)=1$ or $\alpha(x_4)=1$. Therefore, we can take $w=\alpha$.

\medskip

\noindent\textbf{Step 3}. In this step we prove (2). We have
$$\alpha_j(x_1x_2)=\{\frac{ja_1}{r}\}+\{\frac{ja_2}{r}\}\equiv\{\frac{je}{r}\}\mod\mathbb Z,$$
and
$$\alpha_j(x_3x_4)=\{\frac{ja_3}{r}\}+\{\frac{ja_4}{r}\}\equiv\frac{j}{r}\mod\mathbb Z.$$
By (1), there are two cases.

\medskip

\noindent\textbf{Case 1}.  $\alpha_j(x_1x_2)\leq 1$ and $\alpha_j(x_3x_4)>1$. In this case, $\alpha_j(x_1x_2)=\{\frac{je}{r}\}$. By Setting \ref{Setting: before terminal lem}(1.b), $\alpha_j(x_3x_4)<2$, hence $\alpha_j(x_3x_4)=\frac{j}{r}+1$.

\medskip

\noindent\textbf{Case 2}. $\alpha_j(x_1x_2)=\alpha_j(f)+1$ and $\alpha_j(x_3x_4)<1$. In this case, $\alpha_j(x_3x_4)=\frac{j}{r}$. Since $0<\alpha_j(x_1x_2)<2$, $\alpha_j(x_1x_2)=\{\frac{je}{r}\}+1$.

\medskip

\noindent\textbf{Step 4}. In this step we prove (3). By \textbf{Step 1}, $\gcd(a_1,r)=\gcd(a_2,r)=1$. We may assume that $\gcd(e,r)\geq 2$, otherwise (3) follows from Setting \ref{Setting: before terminal lem}(1.a). Let $q:=\frac{r}{\gcd(e,r)}$. 

If $r\mid qa_i$ for some $i\in\{3,4\}$, then $\gcd(e,r)\mid a_i$, and $\gcd(e,r)\mid\gcd(a_i,r)$. Thus $\gcd(a_3,r)=\gcd(e,r)$ or $\gcd(a_4,r)=\gcd(e,r)$ by Setting \ref{Setting: before terminal lem}(1.a), and (3) follows from Setting \ref{Setting: before terminal lem}(1.b). Hence we may assume that  $r\nmid qa_3$ and $r\nmid qa_4$. In particular, $\alpha_q'=\alpha_{r-q}$. There are three cases:

\medskip

\noindent\textbf{Case 1}. $\alpha_q\not\in\Psi$. Then $\alpha_{r-q}\not\in\Psi$. In this case, by (2),
$$\sum_{i=1}^4\{\frac{qa_i}{r}\}=\{\frac{qe}{r}\}+\frac{q}{r}+1=\frac{q}{r}+1$$
and
$$\sum_{i=1}^4\{\frac{(r-q)a_i}{r}\}=\{\frac{(r-q)e}{r}\}+\frac{r-q}{r}+1=\frac{r-q}{r}+1.$$
Thus
$$4=\sum_{i=1}^4(\{\frac{qa_i}{r}\}+\{\frac{(r-q)a_i}{r}\})=3,$$
a contradiction.

\medskip

\noindent\textbf{Case 2}. $\alpha_q\in\Psi_1$. In this case, $\alpha_q=t\beta$ for some $1\leq t\leq k-1$. Since $\beta(x_1x_2)=\beta(f)$ and $\frac{1}{k}<\beta(x_3x_4)\leq\frac{1}{k-1}$, $\alpha_q(x_1x_2)=\alpha_q(f)$ and $$0<\frac{t}{k}<\alpha_q(x_3x_4)\leq\frac{t}{k-1}\leq 1.$$
Since $$\alpha_q(x_3x_4)=\{\frac{qa_3}{r}\}+\{\frac{qa_4}{r}\}\equiv\frac{q}{r}=\frac{1}{\gcd(e,r)}\mod \Zz,$$
we have $\alpha_q(x_3x_4)=\frac{1}{\gcd(e,r)}$, and
$$\frac{t}{k}<\frac{1}{\gcd(e,r)}\leq\frac{t}{k-1}\Longleftrightarrow \frac{k-1}{t}\le\gcd(e,r)<\frac{k}{t}.$$
Thus $\gcd(e,r)=\frac{k-1}{t}$ belongs to a finite set depending only on $k$. Since $q=\frac{r}{\gcd(e,r)}$, $\alpha_q\equiv \frac{1}{\gcd(e,r)}(a_1,a_2,a_3,a_4)\mod\Zz^4$. Thus  $\alpha_q$ belongs to a finite set. Since $\alpha_q=t\beta$ and $1\leq t\leq k-1$, $\beta$ belongs to a finite set, and we are done.

\medskip

\noindent\textbf{Case 3}. $\alpha_q\in\Psi_2$. In this case, $\alpha_{r-q}\in\Psi_1$, hence $\alpha_{r-q}=t\beta$ for some $1\leq t\leq k-1$. Since $\beta(x_1x_2)=\beta(f)$ and $\frac{1}{k}<\beta(x_3x_4)\leq\frac{1}{k-1}$, $\alpha_{r-q}(x_1x_2)=\alpha_{r-q}(f)$ and 
$$0<\frac{t}{k}<\alpha_{r-q}(x_3x_4)\leq\frac{t}{k-1}\leq 1.$$
Since
$$\alpha_{r-q}(x_3x_4)=\{\frac{(r-q)a_3}{r}\}+\{\frac{(r-q)a_4}{r}\}\equiv\frac{r-q}{r}=1-\frac{1}{\gcd(e,r)}\mod \Zz,$$
we have 
$$\frac{k-t}{k}>\frac{1}{\gcd(e,r)}\geq\frac{k-1-t}{k-1},$$
so either $t=k-1$, or $\gcd(e,r)\leq k-1$. There are two sub-cases:

\medskip

\noindent\textbf{Case 3.1}. If $\gcd(e,r)\leq \max\{k-1,6\}$, then $\gcd(e,r)$ belongs to a finite set depending only on $k$. Since $q=\frac{r}{\gcd(e,r)}$, $\alpha_{r-q}\equiv -\frac{1}{\gcd(e,r)}(a_1,a_2,a_3,a_4)\mod\Zz^4$. Thus $\alpha_{r-q}$ belongs to a finite set. Since $\alpha_{r-q}=t\beta$ and $1\leq t\leq k-1$, $\beta$ belongs to a finite set, and we are done.

\medskip

\noindent\textbf{Case 3.2}. If $\gcd(e,r)>\max\{k-1,6\}$, then $t=k-1$, and $q<2q<5q<r$. There are three sub-cases:

\medskip

\noindent\textbf{Case 3.2.1}. There exists $j\in\{2,3,5\}$ such that $\alpha_{jq}\in\Psi_1$. Suppose that $\alpha_{jq}=s\beta$ for some $1\leq s\leq k-1$. Since $\alpha_{r-q}=(k-1)\beta$, $(j(k-1)+s)\beta\equiv\bm{0}\mod \Zz^4$, so $\beta$ belongs to a finite set, and we are done.

\medskip

\noindent\textbf{Case 3.2.2}. There exists $j\in\{2,3,5\}$ such that $\alpha_{jq}\in\Psi_2$. Suppose that $\alpha_{jq}=(s\beta)'$ for some $1\leq s\leq k-1$. Since $\alpha_{r-q}=(k-1)\beta$, $(j(k-1)-s)\beta\equiv\bm{0}\mod \Zz^4$, so 
\begin{itemize}
    \item either $\beta$ belongs to a finite set, in which case we are done, or
    \item $s\equiv j(k-1)\mod r$. In this case, since $1\leq s\leq k-1$, $r$ belongs to a finite set, hence $\beta$ belongs to a finite set, and we are done.
\end{itemize}

\medskip

\noindent\textbf{Case 3.2.3}. For any $j\in\{2,3,5\}$, $\alpha_{jq}\not\in\Psi$. By (2),
$$\sum_{i=1}^4\{\frac{jqa_i}{r}\}=\{\frac{jqe}{r}\}+\frac{jq}{r}+1=\frac{jq}{r}+1$$
and
$$\sum_{i=1}^4\{\frac{(r-jq)a_i}{r}\}=\{\frac{(r-jq)e}{r}\}+\frac{r-jq}{r}+1=\frac{r-jq}{r}+1$$
for any $j\in\{2,3,5\}$, hence
$$2+\sum_{i=3}^4(\{\frac{jqa_i}{r}\}+\{\frac{(r-jq)a_i}{r}\})=\sum_{i=1}^4(\{\frac{jqa_i}{r}\}+\{\frac{(r-jq)a_i}{r}\})=3$$
for any $j\in\{2,3,5\}$. Possibly switching $x_3$ and $x_4$, we may assume that there exist $j_1,j_2\in\{2,3,5\}$ such that $j_1\not=j_2$, $r\mid j_1qa_3$, and $r\mid j_2qa_3$. Thus $r\mid qa_3$, a contradiction.

\medskip

\noindent\textbf{Step 5}. In this step we prove (4). Suppose that $\beta\in N^0$. If $\beta(x_1x_2x_3x_4)-\beta(f)=1$, then since $\sum_{i=1}^4a_i-e\equiv 1\mod r$, we have $\beta\in\mathbb Z^4$, a contradiction. By Setting \ref{Setting: before terminal lem}(3.b), we may assume that $\beta(x_1x_2x_3x_4)-\beta(f)\leq\min\{\frac{12}{13},\frac{1}{k-1}\}$.

By \textbf{Step 1}, we have $\beta(x_1x_2)=\beta(f)$ and $\frac{1}{k}<\beta(x_3x_4)\leq\min\{\frac{12}{13},\frac{1}{k-1}\}$. Thus $\beta(x_3)\not=1$ and $\beta(x_4)\not=1$. By \textbf{Step 1}, $\beta(x_1)\not=1$ and $\beta(x_2)\not=1$. Hence $\beta=\alpha_{k_0}$ for some $1\leq k_0\leq r-1$, and we get (4.a).

By \textbf{Step 1}, $\beta(x_3x_4)\equiv\frac{k_0a_3}{r}+\frac{k_0a_4}{r}\equiv\frac{k_0}{r}\mod\Zz$. Since $\beta(x_3x_4)=\beta(x_1x_2x_3x_4)-\beta(f)\in (\frac{1}{k},\min\{\frac{12}{13},\frac{1}{k-1}\}]$, $\beta(x_3x_4)=\frac{k_0}{r}\in (\frac{1}{k},\min\{\frac{12}{13},\frac{1}{k-1}\}]$. If $\frac{k_0}{r}=\frac{1}{k-1}$, then $\beta=\alpha_{k_0}$ belongs to the finite set $(\frac{1}{k-1}\mathbb Z)^4\cap [0,1]^4$, and we are done. Thus we may assume that $\frac{k_0}{r}<\frac{1}{k-1}$, and we get (4.b).

Now We prove (4.c). Suppose that $\beta(x_1x_2)<1$. Then since $a_1+a_2\equiv e\mod r$, $\{\frac{k_0a_1}{r}\}+\{\frac{k_0a_2}{r}\}=\{\frac{k_0e}{r}\}$. By (4.b), 
$\{\frac{k_0a_3}{r}\}+\{\frac{k_0a_4}{r}\}=\frac{k_0}{r}$. Therefore,
$$\sum_{i=1}^4\{\frac{k_0a_i}{r}\}=\{\frac{k_0e}{r}\}+\frac{k_0}{r}.$$
For any $1\leq j\leq r-1$ such that $j\not=k_0$, there are three cases:

\medskip

\noindent\textbf{Case 1}. $\alpha_j\not\in\Psi$. In this case, by (2) and (4.b), 
$$\sum_{i=1}^4\{\frac{ja_i}{r}\}=\{\frac{je}{r}\}+\frac{j}{r}+1>\{\frac{je}{r}\}+\frac{k_0}{r}+\frac{1}{14}.$$

\medskip

\noindent\textbf{Case 2}. $\alpha_j=t\beta$ for some $2\leq t\leq k-1$. Note that by (4.b), $\frac{k_0}{r}>\frac{1}{k}$. Thus
$$\sum_{i=1}^4\{\frac{ja_i}{r}\}=\sum_{i=1}^4t\cdot\{\frac{k_0a_i}{r}\}=t\cdot\{\frac{k_0e}{r}\}+\frac{tk_0}{r}>\{\frac{tk_0e}{r}\}+\frac{k_0}{r}+\frac{1}{k}=\{\frac{je}{r}\}+\frac{k_0}{r}+\frac{1}{k}.$$

\medskip

\noindent\textbf{Case 3}. $\alpha_j=(t\beta)'$ for some $1\leq t\leq k-1$. In this case, $tk_0\equiv r-j\mod r$. Since $\{\frac{k_0a_3}{r}\}+\{\frac{k_0a_4}{r}\}=\frac{k_0}{r}<\frac{1}{k-1}$, $tk_0=r-j$, and $$\{\frac{tk_0a_3}{r}\}+\{\frac{tk_0a_4}{r}\}\le t\cdot\{\frac{k_0a_3}{r}\}+t\cdot\{\frac{k_0a_4}{r}\}<\frac{t}{k-1}\leq 1.$$
By \textbf{Step 1}, $a_3+a_4\equiv 1\mod r$, so 
$$\{\frac{tk_0a_3}{r}\}+\{\frac{tk_0a_4}{r}\}=\frac{tk_0}{r}.$$ Since $a_1+a_2\equiv e\mod r$, by (3),
$$\{\frac{tk_0a_1}{r}\}+\{\frac{tk_0a_2}{r}\}\leq \{\frac{tk_0e}{r}\}+1.$$
Thus 
$$\sum_{i=1}^4\{\frac{(r-j)a_i}{r}\}=\sum_{i=1}^4\{\frac{tk_0a_i}{r}\}\leq\{\frac{tk_0e}{r}\}+\frac{tk_0}{r}+1.$$
By (3), $\gcd(a_4,r)=\gcd(e,r)$, so $\{\frac{(r-j)a_4}{r}\}+\{\frac{ja_4}{r}\}=\{\frac{(r-j)e}{r}\}+\{\frac{je}{r}\}$. By (4.b), $\frac{k_0}{r}\le \frac{13}{14}$. Thus,
\begin{align*}
\sum_{i=1}^4\{\frac{ja_i}{r}\}=&\sum_{i=1}^4\left(\{\frac{ja_i}{r}\}+\{\frac{(r-j)a_i}{r}\}\right)-\sum_{i=1}^4 \{\frac{(r-j)a_i}{r}\}\\
\geq& 3+\{\frac{(r-j)e}{r}\}+\{\frac{je}{r}\}-\left(\{\frac{tk_0e}{r}\}+\frac{tk_0}{r}+1\right)\\
=&\{\frac{je}{r}\}+\frac{j}{r}+1\\
>&\{\frac{je}{r}\}+\frac{k_0}{r}+\frac{1}{14}.
\end{align*}

By Lemma \ref{lem: transfer to fivefold lemma refined}, $\frac{k_0}{r}$ belongs to a finite set. Since $\beta=\alpha_{k_0}$, $\beta$ belongs to a finite set, and we are done. Therefore, we may assume that $\beta(x_1x_2)\geq 1$. By (4.b), $\{\frac{k_0a_3}{r}\}+\{\frac{k_0a_4}{r}\}=\frac{k_0}{r}$. Since $\beta(x_1x_2)<2$ and $a_1+a_2\equiv e\mod r$, $\{\frac{k_0a_1}{r}\}+\{\frac{k_0a_2}{r}\}=\{\frac{k_0e}{r}\}+1$, and we get (4.c).

\medskip

\noindent\textbf{Step 6}. Finally, we prove (5). If $\beta\not\in N^0$, then $\beta \mod\Zz^4\in \{0,1\}^4$, and $\alpha_j\not\in\Psi$ for any $1\leq j\leq r-1$. In this case, (5) follows from (2). Thus we may assume that $\beta\in N^0$. By (4.a), we may assume that $\beta=\alpha_{k_0}$ for some $1\leq k_0\leq r-1$. If $t\beta\in N^0$ for some $t\geq 2$, then $2\beta\in N^0$, so $2\{\frac{k_0a_1}{r}\}=\{\frac{2k_0a_1}{r}\}$ or $1$, and $2\{\frac{k_0a_2}{r}\}=\{\frac{2k_0a_2}{r}\}$ or $1$. Thus $\{\frac{k_0a_1}{r}\}\leq\frac{1}{2}$ and $\{\frac{k_0a_2}{r}\}\leq\frac{1}{2}$, hence $\{\frac{k_0a_1}{r}\}+\{\frac{k_0a_2}{r}\}\leq 1$. By (4.c), $\{\frac{k_0a_1}{r}\}+\{\frac{k_0a_2}{r}\}=1$, hence $\{\frac{k_0a_1}{r}\}=\{\frac{k_0a_2}{r}\}=\frac{1}{2}$. By (3), $\frac{k_0}{r}=\frac{1}{2}$, so $\beta$ belongs to a finite set, and we are done. Thus we may assume that $t\beta\not\in N^0$ for any $t\geq 2$. Therefore, for any $1\leq j\leq r-1$, there are three cases:

\medskip

\noindent\textbf{Case 1}. $\alpha_j\not\in\Psi$. The equality follows from (2).

\medskip

\noindent\textbf{Case 2}. $\alpha_j=\beta$. The equality follows from (4.c).

\medskip

\noindent\textbf{Case 3}. $\alpha_j=\beta'$. Then $1-\{\frac{k_0a_4}{r}\}=\{\frac{ja_4}{r}\}$. Moreover, since $\beta=\alpha_{k_0}$, $j=r-k_0$. Thus $r\nmid ja_4$, hence $r\nmid je$. By \textbf{Case 2}, we have
$$\sum_{i=1}^4\{\frac{ja_i}{r}\}=4-\sum_{i=1}^4\{\frac{(r-j)a_i}{r}\}=4-(\{\frac{(r-j)e}{r}\}+\frac{r-j}{r}+1)=\{\frac{je}{r}\}+\frac{j}{r}+1.$$

We get (5) and the proof is concluded.
\end{proof}



\begin{prop}\label{prop: cA C case}
Notations and conditions as in Setting \ref{Setting: before terminal lem}. For each positive integer $k$, there exists a finite set $\Ii_k'$ depending only on $k$ satisfying the following. Suppose that $f$ is of cA type, and $\frac{1}{r}(a_1,a_2,a_3,a_4,e)\equiv\frac{1}{r}(a,1,-a,a+1,a+1)\mod \Zz^5$ for some positive integer $a$ such that $\gcd(a,r)=\gcd(a+1,r)=1$. Then either $r\in\Ii_k'$ or $\bm{0}\not=\beta\in\Ii_k'$.
\end{prop}
\begin{proof}
We may assume that $\beta\not\in\Ii_k$ where $\Ii_k$ is the set as in Theorem \ref{thm: cA case up to terminal lemma}.

Since $\gcd(a,r)=\gcd(a+1,r)=1$, we may take an integer $1<b<r$ such that $b(a+1)\equiv 1\mod r$. Then $\gcd(b-1,r)=\gcd(b,r)=1$. We consider
$$\alpha_{b}=\frac{1}{r}(r-b+1,b,b-1,1).$$
 There are three cases.

\medskip

\noindent\textbf{Case 1}. $\alpha\not\in\Psi$. In this case, since $\alpha_b(x_3x_4)<1$, by Theorem \ref{thm: cA case up to terminal lemma}(1), $\alpha_b(f)=\alpha_b(x_1x_2)-1=\frac{1}{r}$. Since $f$ is of cA type, there exists a monomial $\bm{x}\in\mm^2$ such that $\alpha_b(\bm{x})=\alpha_b(f)=\frac{1}{r}$, which is impossible as $\alpha_b(x_3)\geq\alpha_b(x_4)\geq\frac{1}{r}$.

\medskip

\noindent\textbf{Case 2}. $\alpha\in\Psi_1$. Since $\alpha_b(x_4)=\frac{1}{r}$, $\alpha=\beta$. Since $\beta(x_1)+\beta(x_3)=1$ and $\gcd(b-1,r)=1$, either $r=2$ in which case we are done, or $\Psi\cap N^0=\{\beta,\beta'\}$. By Theorem \ref{thm: cA case up to terminal lemma}(4), we may assume that $\frac{1}{k}<\frac{b}{r}<\frac{1}{k-1}$. There are two sub-cases:

\medskip

\noindent\textbf{Case 2.1}. $k=2$. Then $b<r<2b$. There are two sub-cases:

\medskip

\noindent\textbf{Case 2.1.1}. $2b-3\leq r$. Since $\gcd(b-1,r)=1$ and $r<2b$, $r=2b-1$ (resp. $2b-3$), hence 
$$\alpha_{b}=\frac{1}{2b-1}(b,b,b-1,1)\text{ (resp.  }\frac{1}{2b-3}(b-2,b,b-1,1)\text{ )}.$$
We may assume that $b\geq 7$, otherwise $r$ belongs to a finite set and we are done. Then $\alpha_{b+1}=\frac{1}{2b-1}(b+1,b+1,b-2,3)$ (resp. $\alpha_{b+3}=\frac{1}{2b-3}(b-3,b+3,b,3)$).

Since $\beta(x_4)=\alpha_b(x_4)=\frac{1}{r}<1$ and $b\ge7$,  $\alpha_{b+1}\not\in\{\beta,\beta'\}=\Psi\cap N^0$ (resp. $\alpha_{b+3}\not\in\{\beta,\beta'\}=\Psi\cap N^0$). Moreover, $\alpha_{b+1}(x_3x_4)=\frac{b+1}{2b-1}<1$ (resp. $\alpha_{b+3}(x_3x_4)=\frac{b+3}{2b-3}<1$). By Theorem \ref{thm: cA case up to terminal lemma}(1), $\alpha_{b+1}(f)=\alpha_{b+1}(x_1x_2)-1=\frac{3}{2b-1}$ (resp. $\alpha_{b+3}(f)=\alpha_{b+3}(x_1x_2)-1=\frac{3}{2b-3}$).

Thus there exists a monomial $\bm{x}\in\mm^2$ such that $\alpha_{b+1}(\bm{x})=\alpha_{b+1}(f)=\frac{3}{r}$ (resp. $\alpha_{b+3}(\bm{x})=\alpha_{b+3}(f)=\frac{3}{r}$). Since $b\geq 7$, this is impossible as $\alpha_{b+1}(x_3)=\frac{b-2}{2b-1}\geq\frac{3}{2b-1}=\alpha_{b+1}(x_4)$ (resp. $\alpha_{b+3}(x_3)=\frac{b}{2b-3}\geq\frac{3}{2b-3}=\alpha_{b+3}(x_4)$). 

\medskip

\noindent\textbf{Case 2.1.2}. $2b-3>r$. In this case, we consider
$$\alpha_{2b-r}=\frac{1}{r}(2r-2b+2,2b-r,2b-2-r,2).$$
We may assume that $r\geq 7$. Since $\alpha_{2b-r}(x_4)=\frac{2}{r}$ and $\beta(x_4)=\frac{1}{r}$, $\alpha_{2b-r}\not\in\Psi$. Since $b<r$, $\alpha_{2b-r}(x_3x_4)=\frac{2b-r}{r}<1$. By Theorem \ref{thm: cA case up to terminal lemma}(1), $\alpha_{2b-r}(f)=\alpha_{2b-r}(x_1x_2)-1=\frac{2}{r}$. 

Thus there exists a monomial $\bm{x}\in\mm^2$ such that $\alpha_{2b-r}(\bm{x})=\alpha_{2b-r}(f)=\frac{2}{r}$. Since $2b-3>r$, this is impossible as $\alpha_{2b-r}(x_3)=\frac{2b-2-r}{r}\geq\frac{2}{r}=\alpha_{2b-r}(x_4)$.

\medskip

\noindent\textbf{Case 2.2}. $k\geq 3$. Then $kb>r>(k-1)b$. In particular, we may assume that $b>\frac{3}{2}(k+1)$ and $r>k+2$, otherwise $r$ belongs to a finite set and we are done. Let $1\leq c\leq r-1$ be the unique integer such that $c\equiv (k+1)b\mod r$. 

Since $k\ge 3$, $kb>r$ and $b>\frac{3}{2}(k+1)$, we have 
$$2r>2(k-1)b\geq (k+1)b>(k+1)(b-1)>kb.$$ Thus
$$\alpha_c=\frac{1}{r}\big(2r-(k+1)(b-1),(k+1)b-r,(k+1)(b-1)-r,k+1\big).$$
Since $\alpha_c(x_4)+\beta(x_4)=\frac{k+2}{r}$ and $r>k+2$, $\alpha_c\not=\beta'$. Since $k+1\not=1$, $\alpha_c\not=\beta$. Thus $\alpha_c\not\in\Psi$. Since $\alpha_c(x_3x_4)=\frac{(k+1)b-r}{r}<1$, by Theorem \ref{thm: cA case up to terminal lemma}(1), $\alpha_c(f)=\alpha_c(x_1x_2)-1=\frac{k+1}{r}$.

Thus there exists a monomial $\bm{x}\in\mm^2$ such that $\alpha_{c}(\bm{x})=\alpha_{c}(f)=\frac{k+1}{r}$. But this is impossible as $\alpha_c(x_4)=\frac{k+1}{r}$ and $2\alpha_c(x_3)=\frac{2((k+1)(b-1)-r)}{r}>\frac{k+1}{r}$. Note that the last inequality follows from the conditions that $kb>r$ and $b>\frac{3}{2}(k+1)$.

\medskip

\noindent\textbf{Case 3}. $\alpha_b\in\Psi_2$. Since $\gcd(b,r)=\gcd(b-1,r)=1$, $\alpha_b=\beta'$. Thus $\beta=\alpha_{r-b}=\frac{1}{r}(b-1,r-b,r-b+1,r-1)\in N^0$. By Theorem \ref{thm: cA case up to terminal lemma}(4.b), $\frac{2r-b}{r}=\beta(x_3x_4)<\frac{13}{14}<1$, a contradiction.
\end{proof}

\begin{prop}\label{prop: cA D case}
Notations and conditions as in Setting \ref{Setting: before terminal lem}. For each positive integer $k$, there exists a finite set $\Ii_k'$ depending only on $k$ satisfying the following. Suppose that $f$ is of cA type, and $\frac{1}{r}(a_1,a_2,a_3,a_4,e)\equiv\frac{1}{r}(a,-a-1,-a,a+1,-1)\mod \Zz^5$ for some positive integer $a$ such that $\gcd(a,r)=\gcd(a+1,r)=1$. Then either $r\in\Ii_k'$ or $\bm{0}\not=\beta\in\Ii_k'$.
\end{prop}
\begin{proof}
We may assume that $\beta\not\in\Ii_k$ where $\Ii_k$ is the set as in Theorem \ref{thm: cA case up to terminal lemma}.

We consider
$$\alpha_{r-1}=\frac{1}{r}(r-a,a+1,a,r-a-1).$$
There are three cases.

\medskip

\noindent\textbf{Case 1}. $\alpha_{r-1}\not\in\Psi$. Since $\alpha_{r-1}(x_3x_4)=\frac{r-1}{r}<1$, by Theorem \ref{thm: cA case up to terminal lemma}(1), we may assume that $\alpha_{r-1}(f)=\alpha_{r-1}(x_1x_2)-1=\frac{1}{r}$. Thus there exists a monomial $\bm{x}\in\mm^2$ such that $\alpha_{r-1}(\bm{x})=\alpha_{r-1}(f)=\frac{1}{r}$. But this is impossible as $2\alpha_{r-1}(x_3)\geq\frac{1}{r}$ and $\alpha_{r-1}(x_4)\geq\frac{1}{r}$.

\medskip

\noindent\textbf{Case 2}. $\alpha_{r-1}\in\Psi_2$. Since $\gcd(a+1,r)=\gcd(a,r)=1$, $\alpha_{r-1}=\beta'$. Then $\beta=\alpha_1=\frac{1}{r}(a,r-a-1,r-a,a+1)$. By  Theorem \ref{thm: cA case up to terminal lemma}(4), we may assume that $\frac{1}{r}=\beta(x_3x_4)=\frac{r+1}{r}$, which is impossible.

\medskip

\noindent\textbf{Case 3}. $\alpha_{r-1}\in\Psi_1$. Since $\gcd(a+1,r)=\gcd(a,r)=1$, $\alpha_{r-1}=\beta$. By Theorem \ref{thm: cA case up to terminal lemma}(4), we may assume that $\frac{r-1}{r}=\beta(x_3x_4)<\frac{13}{14}$. Thus $r\leq 13$ and we are done.
\end{proof}

\begin{prop}\label{prop: cA B case}
Notations and conditions as in Setting \ref{Setting: before terminal lem}. For each positive integer $k$, there exists a finite set $\Ii_k'$ depending only on $k$ satisfying the following. Suppose that $f$ is of cA type, and $\frac{1}{r}(a_1,a_2,a_3,a_4,e)\equiv\frac{1}{r}(1,a,-a,a+1,a+1)\mod \Zz^5$ for some positive integer $a$ such that $\gcd(a,r)=1$. Then either $r\in\Ii_k'$ or $\bm{0}\not=\beta\in\Ii_k'$. 
\end{prop}
\begin{proof}

We may assume that $\beta\not\in\Ii_k$ where $\Ii_k$ is the set as in Theorem \ref{thm: cA case up to terminal lemma}.

By the ``moreover part'' of Setting \ref{Setting: before terminal lem}, $a+1\not\equiv 0\mod r$. We may assume that $1\leq a<a+1<r$. Moreover, we may assume that $r\geq 5$ otherwise there is nothing left to prove. We consider
$$\alpha_{r-1}:=\frac{1}{r}(r-1,r-a,a,r-a-1).$$
There are three cases.

\medskip

\noindent\textbf{Case 1}. $\alpha_{r-1}\in\Psi_2$. In this case, $\alpha_1=\frac{1}{r}(1,a,r-a,a+1)\in\Psi_1$. Thus $\alpha_1=\beta$. By Theorem \ref{thm: cA case up to terminal lemma}(4.b), $\frac{1}{r}=\alpha_1(x_3x_4)=\frac{r+1}{r}$, which is impossible.

\medskip

\noindent\textbf{Case 2}. $\alpha_{r-1}\not\in\Psi$. In this case, since $\alpha_{r-1}(x_3x_4)<1$, by Theorem \ref{thm: cA case up to terminal lemma}(1), $\alpha_{r-1}(f)=\alpha_{r-1}(x_1x_2)-1=\frac{r-a-1}{r}$. Recall that $f=x_1x_2+g$, and $g\in \mm^2$. There exists a monomial $\bm{x}\in g$ such that $\alpha_{r-1}(\bm{x})=\frac{r-a-1}{r}$. Since $\alpha_{r-1}(x_4)=\frac{r-a-1}{r}$, $\bm{x}=x_3^l$ for some positive integer $l\geq 2$. Since $\alpha_{r-1}(x_3)=\frac{a}{r}$, $l=\frac{r-a-1}{a}$. Thus $r=(l+1)a+1\geq 3a+1\geq 2a+2$. Thus
$$\alpha_{r-2}=\frac{1}{r}(r-2,r-2a,2a,r-2a-2).$$
If $r=2a+2$, then $a=1$ and $r=4$, a contradiction. Thus $r>2a+2$. There are three sub-cases.

\medskip

\noindent\textbf{Case 2.1}. $\alpha_{r-2}\in\Psi_2$. Then $\alpha_2=\frac{1}{r}(2,2a,r-2a,2a+2)\in\Psi_1$. Thus $\alpha_2=\beta$ or $2\beta$. If $\alpha_2=2\beta$, then $\beta=\frac{1}{r}(1,a,\frac{r}{2}-a,a+1)$. By Theorem \ref{thm: cA case up to terminal lemma}(4.a), $\frac{1}{r}(1,a,r-a,a+1)=\alpha_1=\beta$, a contradiction. Thus $\alpha_2=\beta$. By Theorem \ref{thm: cA case up to terminal lemma}(4.b), $\frac{2}{r}=\alpha_2(x_3x_4)=\frac{r+2}{r}$, which is impossible.

\medskip

\noindent\textbf{Case 2.2}. $\alpha_{r-2}\not\in\Psi$. In this case, since $\alpha_{r-2}(x_3x_4)=\frac{r-2}{r}<1$, by Theorem \ref{thm: cA case up to terminal lemma}(1), $\alpha_{r-2}(f)=\alpha_{r-2}(x_1x_2)-1=\frac{r-2a-2}{r}$. Thus there exists a monomial $\bm{y}\in g$ such that $\alpha_{r-2}(\bm{y})=\alpha_{r-2}(f)=\frac{r-2a-2}{r}$. Since $\alpha_{r-2}(x_4)=\frac{r-2a-2}{r}$. Thus $\bm{y}=x_3^s$ for some integer $s\geq 2$. Since $x_3^l\in g$, we have 
$$\alpha_{r-2}(x_3^s)=\alpha_{r-2}(f)\leq\alpha_{r-2}(x_3^l),$$
and
$$\alpha_{r-1}(x_3^s)\geq\alpha_{r-1}(f)=\alpha_{r-1}(x_3^l),$$
so $l=s$. Thus $l=\frac{r-2a-2}{2a}=\frac{r-a-1}{a}$, a contradiction.

\medskip

\noindent\textbf{Case 2.3}. $\alpha_{r-2}\in\Psi_1$. Since $\gcd(a,r)=1$, $\alpha_{r-2}=\beta$ or $2\beta$. There are two sub-cases.

\medskip

\noindent\textbf{Case 2.3.1}. $\alpha_{r-2}=\beta$. In this sub-case, by Theorem \ref{thm: cA case up to terminal lemma}(4.b), $\frac{r-2}{r}=\beta(x_3x_4)<\frac{13}{14}$, so $r\leq 27$ and we are done.

\medskip

\noindent\textbf{Case 2.3.2}. $\alpha_{r-2}=2\beta$. In this sub-case, by Theorem \ref{thm: cA case up to terminal lemma}(4.a),
$$\beta=\alpha_{\frac{r}{2}-1}=\frac{1}{r}(\frac{r}{2}-1,\frac{r}{2}-a,a,\frac{r}{2}-a-1).$$
Thus $\beta(x_1x_2)=\frac{r-a-1}{r}<1$, which contradicts Theorem \ref{thm: cA case up to terminal lemma}(4.c).

\medskip

\noindent\textbf{Case 3}. $\alpha_{r-1}\in\Psi_1$. In this case, $\alpha_{r-1}=\beta$ as $r-a-(r-a-1)=1$. By Theorem \ref{thm: cA case up to terminal lemma}(4.b), $\beta(x_3x_4)=\frac{r-1}{r}<\frac{13}{14}$. Hence $r\leq 13$, and we are done.
\end{proof}

\begin{proof}[Proof of Theorem \ref{thm: cA type beta finite}] By Theorem \ref{thm: cA case up to terminal lemma}(3,5), $a_i,e$ satisfy the conditions of terminal lemma. Since $xy\in f$, $a_1+a_2\equiv e \mod r$. Recall that by the ``moreover part'' of Setting \ref{Setting: before terminal lem},  $\frac{1}{r}(a_1,a_2,a_3,a_4,e)\not\equiv\frac{1}{r}(a,-a,1,0,0)\mod \Zz^5$ for any for any integer $a$ such that $\gcd(a,r)=1$. Thus by the terminal lemma (Theorem \ref{thm: terminal lemma}),  possibly interchanging $a_1,a_2$ or $a_3,a_4$, one of the following holds.
\begin{enumerate}
    \item $\frac{1}{r}(a_1,a_2,a_3,a_4,e)\equiv\frac{1}{r}(a,1,-a,a+1,a+1) \mod \Zz^5$ for some positive integer $a$ such that $\gcd(a,r)=\gcd(a+1,r)=1$.
    \item $\frac{1}{r}(a_1,a_2,a_3,a_4,e)\equiv\frac{1}{r}(a,-a-1,-a,a+1,-1) \mod \Zz^5$ for some positive integer $a$ such that $\gcd(a,r)=\gcd(a+1,r)=1$.
    
    \item $\frac{1}{r}(a_1,a_2,a_3,a_4,e)\equiv\frac{1}{r}(1,a,-a,a+1,a+1) \mod \Zz^5$ for some integer $a$ such that $\gcd(a,r)=1$ and $\gcd(a+1,r)>1$.
\end{enumerate}
Now the theorem follows from Propositions \ref{prop: cA C case}, \ref{prop: cA D case}, and \ref{prop: cA B case}. 
\end{proof}

\section{Non-cA type}\label{note: noncA type}
The goal of this subsection is to show Theorem \ref{thm: non cA type beta finite}.

\begin{thm}\label{thm: non-cA case up to terminal lemma}
Notations and conditions as in Setting \ref{Setting: before terminal lem}. Suppose that $f$ is not of cA type (Setting \ref{Setting: before terminal lem}(2.a)). Then there exists a finite set $\Ii_k\subset\mathbb Q^4\cap [0,1]$ depending only on $k$ satisfying the following. Then possibly switching $x_2,x_3$ and $x_4$, either $\bm{0}\not=\beta\in\Ii_k$, or we have the following.
\begin{enumerate}
   \item For any $\alpha\in N^0\backslash\Psi$, there exists $w\in\{\alpha,\alpha'\}$, such that
    \begin{enumerate}
        \item $w(f)=2w(x_1)\leq 1$ and $w'(f)=2w'(x_1)-1\geq 0$,
        \item $w(x_2x_3x_4)>1+w(x_1)$ and $w'(x_2x_3x_4)<1+w'(x_1)$, and
        \item $2w(x_1)=1$ if and only if $2w'(x_1)-1=0$. Moreover, if $2w(x_1)=1$, then $w'\equiv a_{\frac{r}{2}}\mod \Zz^4$, and there exists $\{i_2,i_3,i_4\}=\{2,3,4\}$ such that $(w(x_{i_2}),w(x_{i_3}),w(x_{i_4}))=(\frac{1}{2},\frac{1}{2},1)$ and $(w'(x_{i_2}),w'(x_{i_3}),w'(x_{i_4}))=(\frac{1}{2},\frac{1}{2},0)$.
    \end{enumerate}
    \item For any $1\leq j\leq r-1$ such that $\alpha_j\not\in\Psi$, either
    \begin{itemize}
        \item  $2\{\frac{ja_1}{r}\}=\{\frac{je}{r}\}$ and $\{\frac{ja_2}{r}\}+\{\frac{ja_3}{r}\}+\{\frac{ja_4}{r}\}=\{\frac{ja_1}{r}\}+\frac{j}{r}+1$, or
        \item $2\{\frac{ja_1}{r}\}=\{\frac{je}{r}\}+1$ and $\{\frac{ja_2}{r}\}+\{\frac{ja_3}{r}\}+\{\frac{ja_4}{r}\}=\{\frac{ja_1}{r}\}+\frac{j}{r}$.
    \end{itemize}
    \item One of the following holds:
    \begin{enumerate}
        \item $\gcd(a_1,r)=\gcd(e,r)\geq 2$ and $\gcd(a_2,r)=\gcd(a_3,r)=\gcd(a_4,r)=1$.
        \item $2\nmid r$ and $\gcd(a_i,r)=\gcd(e,r)=1$ for any $1\leq i\leq 4$.
        \item $\gcd(a_4,r)=\gcd(e,r)=2$, and $\gcd(a_1,r)=\gcd(a_2,r)=\gcd(a_3,r)=1$.
    \end{enumerate}
    \item If $\beta\in N^0$, then there exists $1\leq k_0\leq r-1$, such that 
    \begin{enumerate}
        \item $\beta\equiv\alpha_{k_0}\mod \Zz^4$,
        \item $\frac{1}{k}<\frac{k_0}{r}<\min\{\frac{13}{14},\frac{1}{k-1}\}$, and
            \item if $\beta=\alpha_{k_0}$, then $\beta(f)=2\beta(x_1)\geq 1$, $2\{\frac{k_0a_1}{r}\}=\{\frac{k_0e}{r}\}+1$, and $\{\frac{k_0a_2}{r}\}+\{\frac{k_0a_3}{r}\}+\{\frac{k_0a_4}{r}\}=\{\frac{k_0a_1}{r}\}+\frac{k_0}{r}$.
    \end{enumerate}
    \item  For any $1\leq j\leq r-1$,
    $$\sum_{i=1}^4\{\frac{ja_i}{r}\}=\{\frac{je}{r}\}+\frac{j}{r}+1.$$
    \item If $\gcd(a_1,r)=\gcd(e,r)\geq 2$, then $\gcd(a_1,r)=\gcd(e,r)=r$.
\end{enumerate}
\end{thm}
\begin{proof}
\noindent\textbf{Step 1}. In this step we summarize some auxiliary results that will be used later.

By Setting \ref{Setting: before terminal lem}(2), $x_1^2\in f$, so $2a_1\equiv e\mod r$, and $a_2+a_3+a_4\equiv a_1+1\mod r$. Thus $\alpha(f)\equiv 2\alpha(x_1)\mod \Zz$ for any $\alpha\in N$. 

For any $\alpha\in N^0$, $0\leq\alpha(f)\leq 2\alpha(x_1)\leq 2$. Thus $\alpha(f)\in\{2\alpha(x_1),2\alpha(x_1)-1,2\alpha(x_1)-2\}$. If $\alpha(f)=2\alpha(x_1)-2$, then $\alpha(x_1)=1$ and $\alpha(f)=0$, hence $\alpha(x_i)=0$ for some $i\in\{2,3,4\}$. Since $\alpha\in N$, there exists $l\in [1,r-1]\cap\Zz_{\ge 1}$ such that $\alpha\equiv\alpha_l\mod \Zz^4$, thus $r\mid (la_1,la_i)$ which contradicts Setting \ref{Setting: before terminal lem}(1.b). Therefore, $\alpha(f)\in\{2\alpha(x_1),2\alpha(x_1)-1\}$. If $\alpha(x_1x_2x_3x_4)-\alpha(f)>1$ and $\alpha(f)=2\alpha(x_1)$, then $\alpha(x_2x_3x_4)>\alpha(x_1)+1$.

If $k\neq 1$, then since $\frac{1}{k}<\beta(x_1x_2x_3x_4)-\beta(f)\leq\frac{1}{k-1}$, either
\begin{itemize}
    \item $\beta(f)=2\beta(x_1)$ and $\frac{1}{k}<\beta(x_2x_3x_4)-\beta(x_1)\leq\frac{1}{k-1}$, or
    \item $\beta(f)=2\beta(x_1)-1$ and $-1+\frac{1}{k}<\beta(x_2x_3x_4)-\beta(x_1)\leq -1+\frac{1}{k-1}$.
\end{itemize}

Finally, since switching $x_2,x_3,$ and $x_4$ will not influence (1)(2)(4)(5), we will only have a possibly switching of $x_2,x_3,x_4$ when we prove (3).

\medskip

\noindent\textbf{Step 2}. In this step we prove (1). Pick $\alpha\in N^0\backslash\Psi$, then $\alpha'\in N^0\backslash\Psi$. By \textbf{Step 1}, there are two cases:

\medskip

\noindent\textbf{Case 1}. $\alpha(f)=2\alpha(x_1)-1$. In this case, $2\alpha(x_1)\geq 1$. There are two sub-cases.

\medskip

\noindent\textbf{Case 1.1}. $2\alpha(x_1)=1$. Then $\alpha(x_1)=\frac{1}{2}$, so $\alpha'(x_1)=\frac{1}{2}$. Moreover, $\alpha(f)=0$, hence there exists $i\in\{2,3,4\}$ such that $\alpha(x_i)=0$. Thus $\alpha'(x_i)=1$. Since $\alpha\in N$, there exists $l\in [1,r-1]\cap\Zz_{\ge 1}$ such that $\alpha\equiv\alpha_l\mod \Zz^4$. Thus $r\mid la_i$, $r\mid 2la_1$, and $r\mid 2l\gcd(a_1,a_i)$. By Setting \ref{Setting: before terminal lem}(1.b), $r,a_i$ are even and $l=\frac{r}{2}$. By Setting \ref{Setting: before terminal lem}(1.b), $2\nmid a_j$ for any $j\not=i$. Since $\alpha\in N^0$, $\alpha(x_j)=\frac{1}{2}$ for any $j\not=i$. Thus $\alpha'(x_j)=\frac{1}{2}$ for any $j\not=i$. Now $\alpha(x_2x_3x_4)=1<\frac{3}{2}=1+\alpha(x_1)$, and $\alpha'(x_2x_3x_4)=2>\frac{3}{2}=1+\alpha'(x_1)$. Thus we may let $w=\alpha'$.

\medskip

\noindent\textbf{Case 1.2}. $2\alpha(x_1)>1$. Then $2\alpha'(x_1)<1$. By \textbf{Step 1}, $\alpha'(f)=2\alpha'(x_1)$ and $\alpha'(x_2x_3x_4)>\alpha'(x_1)+1$. Thus
$$\alpha(x_2x_3x_4)=3-\alpha'(x_2x_3x_4)<2-\alpha'(x_1)=\alpha(x_1)+1.$$
Thus we may take $w=\alpha'$.

\medskip

\noindent\textbf{Case 2}. $\alpha(f)=2\alpha(x_1)$. By \textbf{Step 1}, $\alpha(x_2x_3x_4)>\alpha(x_1)+1$, hence
$$\alpha'(x_2x_3x_4)=3-\alpha(x_2x_3x_4)<2-\alpha(x_1)=\alpha'(x_1)+1.$$
By \textbf{Step 1}, $\alpha'(f)\not=2\alpha'(x_1)$, hence $\alpha'(f)=2\alpha'(x_1)-1$. Thus $\alpha'$ satisfies \textbf{Case 1}. By \textbf{Case 1}, we may take $w=\alpha$.

\medskip

\noindent\textbf{Step 3}. In this step we prove (2). Pick $j\in [1,r-1]\cap\Zz_{\ge 1}$ such that $\alpha_j\not\in\Psi$. By \textbf{Step 1}, $2\frac{ja_1}{r}\equiv\frac{je}{r}\mod \Zz$ and $\frac{ja_2}{r}+\frac{ja_3}{r}+\frac{ja_4}{r}\equiv\frac{ja_1}{r}+\frac{j}{r}\mod \Zz$. By (1), there are two cases.

\medskip

\noindent\textbf{Case 1}. $\alpha_j(f)=2\alpha_j(x_1)\leq 1$ and $\alpha_j(x_2x_3x_4)>1+\alpha_j(x_1)$. There are two sub-cases.

\medskip

\noindent\textbf{Case 1.1}. $2\alpha_j(x_1)=1$. In this case, by (1.c), $j=\frac{r}{2}$. In particular, $je\equiv 2ja_1\equiv 0 \mod r$. Then 
$$2\{\frac{ja_1}{r}\}=\{\frac{je}{r}\}+1,\, \{\frac{ja_2}{r}\}+\{\frac{ja_3}{r}\}+\{\frac{ja_4}{r}\}=\{\frac{ja_1}{r}\}+\frac{j}{r},$$
and (2) follows.

\medskip

\noindent\textbf{Case 1.2}. $2\alpha_j(x_1)<1$. Then $2\{\frac{ja_1}{r}\}=\{\frac{je}{r}\}$. By (1),
$$1+\{\frac{ja_1}{r}\}=1+\alpha_j(x_1)<\alpha_j(x_2x_3x_4)=\{\frac{ja_2}{r}\}+\{\frac{ja_3}{r}\}+\{\frac{ja_4}{r}\}<3,$$
thus $\alpha_j(x_2x_3x_4)\in\{1+\frac{j}{r}+\{\frac{ja_1}{r}\},2+\frac{j}{r}+\{\frac{ja_1}{r}\}\}$. If $\alpha_j(x_2x_3x_4)=1+\frac{j}{r}+\{\frac{ja_1}{r}\}$ then we are done, so we may assume that $\alpha_j(x_2x_3x_4)=2+\frac{j}{r}+\{\frac{ja_1}{r}\}$. Thus $\alpha_j'(x_2x_3x_4)=1-\frac{j}{r}-\{\frac{ja_1}{r}\}$ and $\alpha_j'(x_1)=1-\{\frac{ja_1}{r}\}$. By (1), $\alpha_j'(f)=2\alpha'(x_1)-1=1-2\{\frac{ja_1}{r}\}$. Thus $\alpha_j'(x_1x_2x_3x_4)-\alpha'(f)=1-\frac{j}{r}\leq 1$, a contradiction.

\medskip

\noindent\textbf{Case 2}. $\alpha_j(f)=2\alpha_j(x_1)-1$ and $\alpha_j(x_2x_3x_4)<1+\alpha_j(x_1)$. Then $2>2\alpha_j(x_1)\geq 1$, hence $2\alpha_j(x_1)=1+\{\frac{je}{r}\}$. Moreover, by Setting \ref{Setting: before terminal lem}(3.a) and (1.b), $1+\alpha_j(x_1)>\alpha_j(x_2x_3x_4)>\alpha_j(x_1)$. Thus $\alpha_j(x_2x_3x_4)=\{\frac{ja_2}{r}\}+\{\frac{ja_3}{r}\}+\{\frac{ja_4}{r}\}=\{\frac{ja_1}{r}\}+\frac{j}{r}$, and (2) follows.

\medskip

\noindent\textbf{Step 4}. In this step we prove (3). We first prove the following claim.

\begin{claim}\label{claim: non-cA gcd a e}
If $\gcd(a_1,r)\geq 2$, then either $\beta$ belongs to a finite set depending only on $k$, or $\gcd(a_1,r)=\gcd(e,r)$.
\end{claim}
\begin{proof}
Suppose the claim does not hold, then since $e\equiv 2a_1\mod r$, $\gcd(e,r)=2\gcd(a_1,r)$. In particular, $2\mid r$, $2\nmid \frac{a_1}{\gcd(a_1,r)}$. Let $q:=\frac{r}{\gcd(e,r)}=\frac{r}{2\gcd(a_1,r)}$. Then $r\mid qe$ and $qa_1\equiv (r-q)a_1\equiv \frac{r}{2}\mod r$. By Setting \ref{Setting: before terminal lem}(1.b), $\frac{r}{2}\nmid qa_i$ for any $i\in\{2,3,4\}$. Thus $\alpha_q'=\alpha_{r-q}$. There are three cases.

\medskip

\noindent\textbf{Case 1}. $\alpha_q\not\in\Psi$. In this case, by (2),
$$\sum_{i=1}^4\{\frac{qa_i}{r}\}=\{\frac{qe}{r}\}+\frac{q}{r}+1=\frac{q}{r}+1$$
and
$$\sum_{i=1}^4\{\frac{(r-q)a_i}{r}\}=\{\frac{(r-q)e}{r}\}+\frac{r-q}{r}+1=\frac{r-q}{r}+1,$$
hence
$$4=\sum_{i=1}^4(\{\frac{qa_i}{r}\}+\{\frac{(r-q)a_i}{r}\})=3,$$
a contradiction.

\medskip

\noindent\textbf{Case 2}. $\alpha_q\in\Psi_1$. In this case, $\alpha_q=t\beta$ for some $1\leq t\leq k-1$. Recall that $\alpha_q(f)\equiv \frac{qe}{r}\mod \Zz,$ and $a_1+a_2+a_3+a_4-e\equiv 1 \mod r$, we have 
$$t(\beta(x_1x_2x_3x_4)-\beta(f))=\alpha_q(x_1x_2x_3x_4)-\alpha_q(f)\equiv\frac{q}{r}=\frac{1}{2\gcd(a_1,r)}\mod \Zz.$$
Since
$$t(\beta(x_1x_2x_3x_4)-\beta(f))\in (\frac{t}{k},\frac{t}{k-1}],$$
we have
$$\frac{t}{k}<\frac{1}{2\gcd(a_1,r)}\leq\frac{t}{k-1},$$
so $\gcd(a_1,r)=\frac{k-1}{2t}$ belongs to a finite set depending only on $k$. Since $\alpha_q\equiv\frac{1}{2\gcd(a_1,r)}(a_1,a_2,a_3,a_4)\mod \Zz^4$, $\alpha_q$ belongs to a finite set. Since $\alpha_q=t\beta$ and $1\leq t\leq k-1$, $\beta$ belongs to a finite set, a contradiction.

\medskip

\noindent\textbf{Case 3}. $\alpha_q\in\Psi_2$. In this case, $\alpha_{r-q}\in\Psi_1$, hence $\alpha_{r-q}=t\beta$ for some $1\leq t\leq k-1$. Since $$t(\beta(x_1x_2x_3x_4)-\beta(f))=\alpha_{r-q}(x_1x_2x_3x_4)-\alpha_{r-q}(f)\equiv\frac{r-q}{r}=1-\frac{1}{2\gcd(a_1,r)}\mod \Zz$$
and
$$t(\beta(x_1x_2x_3x_4)-\beta(f))\in (\frac{t}{k},\frac{t}{k-1}],$$
we have that
$$\frac{k-t}{k}>\frac{1}{2\gcd(a_1,r)}\geq\frac{k-1-t}{k-1},$$
so either $t=k-1$, or $2\gcd(a_1,r)\leq k-1$. There are two sub-cases:

\medskip

\noindent\textbf{Case 3.1}. If $2\gcd(a_1,r)\leq\max\{k-1,12\}$, then $\gcd(a_1,r)$ belongs to a finite set depending only on $k$. Since $q=\frac{r}{2\gcd(a_1,r)}$, $\alpha_{r-q}\equiv-\frac{1}{2\gcd(a_1,r)}(a_1,a_2,a_3,a_4)\mod\Zz^4$. Thus $\alpha_{r-q}$ belongs to a finite set. Since $\alpha_{r-q}=t\beta$ and $1\leq t\leq k-1$, $\beta$ belongs to a finite set, a contradiction.

\medskip

\noindent\textbf{Case 3.2}. If $2\gcd(a_1,r)>\max\{k-1,12\}$, then $t=k-1$, and $q<3q<5q<7q<11q<r$. There are three sub-cases:

\medskip

\noindent\textbf{Case 3.2.1}. There exists $j\in\{3,5,7,11\}$ such that $\alpha_{jq}\in\Psi_1$. Suppose that $\alpha_{jq}=s\beta$ for some $1\leq s\leq k-1$. Since $\alpha_{r-q}=(k-1)\beta$, $(j(k-1)+s)\beta\equiv\bm{0}\mod \Zz^4$, so $\beta$ belongs to a finite set, a contradiction. 

\medskip

\noindent\textbf{Case 3.2.2}. There exists $j\in\{3,5,7,11\}$ such that $\alpha_{jq}\in\Psi_2$. Suppose that $\alpha_{jq}=(s\beta)'$ for some $1\leq s\leq k-1$. Since  $\alpha_{r-q}=(k-1)\beta$, $(j(k-1)-s)\beta\equiv\bm{0}\mod \Zz^4$, so either 
\begin{itemize}
    \item $\beta$ belongs to a finite set, in which case we get a contradiction, or
    \item $s\equiv j(k-1)\mod r$. In this case, since $1\leq s\leq k-1$, $r$ belongs to a finite set, hence $\beta$ belongs to a finite set, and we get a contradiction again.
\end{itemize}

\medskip

\noindent\textbf{Case 3.2.3}. For any $j\in\{3,5,7,11\}$, $\alpha_{jq}\not\in\Psi$. By (2),
$$\sum_{i=2}^4\{\frac{jqa_i}{r}\}=\frac{1}{2}+\frac{jq}{r}$$
and
$$\sum_{i=2}^4\{\frac{(r-jq)a_i}{r}\}=\frac{1}{2}+\frac{(r-jq)}{r}$$
for any $j\in\{3,5,7,11\}$, hence
$$\sum_{i=2}^4(\{\frac{jqa_i}{r}\}+\{\frac{(r-jq)a_i}{r}\})=2$$
for any $j\in\{3,5,7,11\}$. Possibly switching $x_2,x_3,x_4$, we may assume that there exist $j_1,j_2\in\{3,5,7,11\}$ such that $j_1\not=j_2$, $r\mid j_1qa_2$, and $r\mid j_2qa_2$. Thus $r\mid qa_2$, a contradiction.
\end{proof}
\noindent\textit{Proof of Theorem \ref{thm: non-cA case up to terminal lemma} continued}. We continue \textbf{Step 4}. 

If $\gcd(a_1,r)\geq 2$, then by Claim \ref{claim: non-cA gcd a e}, we may assume that $\gcd(a_1,r)=\gcd(e,r)$. By Setting \ref{Setting: before terminal lem}(1.a-b), $\gcd(a_i,r)=1$ for $2\le i\le 4$, and (3.a) holds. So we may assume that $\gcd(a_1,r)=1$. Since $e\equiv 2a_1\mod r$, $\gcd(e,r)=1$ or $2$. If $\gcd(e,r)=1$, then by Setting \ref{Setting: before terminal lem}(1.a), (3.b) holds. Thus we may assume that $\gcd(e,r)=2$. There are two cases.

\medskip

\noindent\textbf{Case 1}. $a_2,a_3,a_4$ are odd. In this case, by Setting \ref{Setting: before terminal lem}(1.a), we have $\alpha_{\frac{r}{2}}=(\frac{1}{2},\frac{1}{2},\frac{1}{2},\frac{1}{2})$. By (1.b), $\alpha_{\frac{r}{2}}\in\Psi$. Thus $\beta$ belongs to the finite set $$\{\frac{1}{2t}(1,1,1,1),\frac{1}{2t}(2t-1,2t-1,2t-1,2t-1)\mid 1\leq t\leq k-1\},$$
and we are done. 

\medskip

\noindent\textbf{Case 2}. $2\mid a_i$ for some $i\in\{2,3,4\}$. By Setting \ref{Setting: before terminal lem}(1.b), $2\nmid a_j$ for any $j\not=i$. Possibly switching $x_2,x_3,x_4$, we may assume that $2\mid a_4$. By Setting \ref{Setting: before terminal lem}(1.a), $\gcd(a_1,r)=\gcd(a_2,r)=\gcd(a_3,r)=1,\gcd(a_4,r)=\gcd(e,r)=2$, and (3.c) holds.

\medskip

\noindent\textbf{Step 5}. In this step we prove (4). (4.a) follows from the construction of $\beta$. Suppose that $\beta\in N^0$. If $\beta(x_1x_2x_3x_4)-\beta(f)=1$, then since $\sum_{i=1}^4a_i-e\equiv 1\mod r$, $\beta\in\mathbb Z^4$, a contradiction. By Setting \ref{Setting: before terminal lem}(3.b), we may assume that $\beta(x_1x_2x_3x_4)-\beta(f)\leq\min\{\frac{12}{13},\frac{1}{k-1}\}$.

Since $\sum_{i=1}^4a_i-e\equiv 1\mod r$, $$\frac{1}{k}<\frac{k_0}{r}=\beta(x_1x_2x_3x_4)-\beta(f)\leq\min\{\frac{12}{13},\frac{1}{k-1}\}.$$
To prove (4.b), we only need to deal with the case when $\frac{k_0}{r}=\frac{1}{k-1}$. In this case, $\alpha_{k_0}$ belongs to a finite set, hence $\beta$ belongs to a finite set, and we are done. 

We prove (4.c). Suppose that $\beta=\alpha_{k_0}$. By \textbf{Step 1}, there are two cases.

\medskip

\noindent\textbf{Case 1}. $\beta(f)=2\beta(x_1)-1$. Then $-1+\frac{1}{k}<\beta(x_2x_3x_4)-\beta(x_1)\leq\min\{-\frac{1}{13},-1+\frac{1}{k-1}\}$, and $2\beta(x_1)\geq 1$. Since $2a_1\equiv e\mod r$, $2\{\frac{k_0a_1}{r}\}=\{\frac{k_0e}{r}\}+1$. Since $a_2+a_3+a_4\equiv a_1+1\mod r$, $\{\frac{k_0a_2}{r}\}+\{\frac{k_0a_3}{r}\}+\{\frac{k_0a_4}{r}\}\equiv \{\frac{k_0a_1}{r}\}+\frac{k_0}{r}\mod \Zz$, hence $\{\frac{k_0a_2}{r}\}+\{\frac{k_0a_3}{r}\}+\{\frac{k_0a_4}{r}\}=\{\frac{k_0a_1}{r}\}+\frac{k_0}{r}-1$. Thus
$$\sum_{i=1}^4\{\frac{k_0a_i}{r}\}=\{\frac{k_0e}{r}\}+\frac{k_0}{r}.$$
Since $2\beta(x_1)\geq 1$, either $2\beta(x_1)=1$, or $\Psi\cap N^0=\{\beta,\beta'\}$. If $2\beta(x_1)=1$, then $\beta(f)=0$, hence $\beta(x_i)=0$ for some $i\in\{2,3,4\}$. By (3), $r\leq 2$, hence $\beta$ belongs to a finite set, and we are done. Thus we have $\Psi\cap N^0=\{\beta,\beta'\}$.

For any $1\leq j\leq r-1$ such that $j\not=k_0$, if  $\alpha_j\not\in\Psi$, then by (2),
$$\sum_{i=1}^4\{\frac{ja_i}{r}\}=\{\frac{je}{r}\}+\frac{j}{r}+1\geq\{\frac{je}{r}\}+\frac{j}{r}+\frac{1}{14}.$$
Otherwise, $\alpha_j=\beta'$. Thus $j=r-k_0$. If there exists $i$ such that $r\mid k_0a_i$, then $r\mid k_0e$ and $r\nmid k_0a_j$ for any $j\not=i$ by (3). We have
\begin{align*}
&\sum_{i=1}^4\{\frac{ja_i}{r}\}=\sum_{i=1}^4\{\frac{(r-k_0)a_i}{r}\}\\
=&3-\sum_{i=1}^4\{\frac{k_0a_i}{r}\}=3-(\{\frac{k_0e}{r}\}+\frac{k_0}{r})\\
=&2+\frac{j}{r}>\{\frac{je}{r}\}+\frac{k_0}{r}+1.
\end{align*}
If $r\nmid k_0a_i$ for each $i$, then
\begin{align*}
&\sum_{i=1}^4\{\frac{ja_i}{r}\}=\sum_{i=1}^4\{\frac{(r-k_0)a_i}{r}\}\\
=&4-\sum_{i=1}^4\{\frac{k_0a_i}{r}\}=4-(\{\frac{k_0e}{r}\}+\frac{k_0}{r})\\
\geq &2+\{\frac{je}{r}\}+\frac{j}{r}>\{\frac{je}{r}\}+\frac{k_0}{r}+1.
\end{align*}
By Lemma \ref{lem: transfer to fivefold lemma refined}, $\frac{k_0}{r}$ belongs to a finite set, hence $\beta=\alpha_{k_0}$ belongs to a finite set, and we are done.

\medskip

\noindent\textbf{Case 2}. $\beta(f)=2\beta(x_1)$. Then $\frac{1}{k}<\beta(x_2x_3x_4)-\beta(x_1)\leq\min\{\frac{12}{13},\frac{1}{k-1}\}$. Since $a_2+a_3+a_4\equiv a_1+1\mod r$, $\{\frac{k_0a_2}{r}\}+\{\frac{k_0a_3}{r}\}+\{\frac{k_0a_4}{r}\}\equiv \{\frac{k_0a_1}{r}\}+\frac{k_0}{r}\mod \Zz$, hence $\{\frac{k_0a_2}{r}\}+\{\frac{k_0a_3}{r}\}+\{\frac{k_0a_4}{r}\}=\{\frac{k_0a_1}{r}\}+\frac{k_0}{r}$. In particular, $\beta(x_2x_3x_4)-\beta(x_1)=\frac{k_0}{r}$. To prove (4.c), we only need to show that $\beta$ belongs to a finite set when $2\{\frac{k_0a_1}{r}\}=\{\frac{k_0e}{r}\}$. In this case,
$$\sum_{i=1}^4\{\frac{k_0a_i}{r}\}=\{\frac{k_0e}{r}\}+\frac{k_0}{r}.$$
For any $1\leq j\leq r-1$ such that $j\not=k_0$, there are three sub-cases.

\medskip

\noindent\textbf{Case 2.1}. $\alpha_j\not\in\Psi$. In this case, by (2) and (4.b),
$$\sum_{i=1}^4\{\frac{ja_i}{r}\}=\{\frac{je}{r}\}+\frac{j}{r}+1>\{\frac{je}{r}\}+\frac{k_0}{r}+\frac{1}{14}.$$

\medskip

\noindent\textbf{Case 2.2}. $\alpha_j=t\beta$ for some $2\leq t\leq k-1$. In this case, by (4.b),
$$\sum_{i=1}^4\{\frac{ja_i}{r}\}=t\sum_{i=1}^4\{\frac{k_0a_i}{r}\}=t(\{\frac{k_0e}{r}\}+\frac{k_0}{r})>\{\frac{tk_0e}{r}\}+\frac{k_0}{r}+\frac{1}{k}=\{\frac{je}{r}\}+\frac{k_0}{r}+\frac{1}{k}.$$

\medskip

\noindent\textbf{Case 2.3}. $\alpha_j=(t\beta)'$ for some $1\leq t\leq k-1$. In this case, $tk_0\equiv r-j\mod r$. Since $\frac{1}{k}<\frac{k_0}{r}<\frac{1}{k-1}$, $0<tk_0<r$. Since $1\leq j\leq r-1$, $tk_0=r-j$.

Since $t\{\frac{k_0a_i}{r}\}=\{\frac{tk_0a_i}{r}\}$ for any $i\in\{1,2,3,4\}$, we have
$$\{\frac{tk_0a_2}{r}\}+\{\frac{tk_0a_3}{r}\}+\{\frac{tk_0a_4}{r}\}=\{\frac{tk_0a_1}{r}\}+\frac{tk_0}{r}.$$
Thus by (3),
\begin{align*}
&\{\frac{ja_2}{r}\}+\{\frac{ja_3}{r}\}+\{\frac{ja_4}{r}\}\\
=&3-(\{\frac{tk_0a_2}{r}\}+\{\frac{tk_0a_3}{r}\}+\{\frac{tk_0a_4}{r}\})\\
=&3-(\{\frac{tk_0a_1}{r}\}+\frac{tk_0}{r})\geq 1+\{\frac{ja_1}{r}\}+\frac{j}{r}.
\end{align*}
Since $2a_1\equiv e\mod r$, $2\{\frac{ja_1}{r}\}\geq\{\frac{je}{r}\}$. Thus
$$\sum_{i=1}^4\{\frac{ja_i}{r}\}\geq\{\frac{je}{r}\}+\frac{j}{r}+1>\{\frac{je}{r}\}+\frac{k_0}{r}+\frac{1}{14}.$$
By Lemma \ref{lem: transfer to fivefold lemma refined}, $\frac{k_0}{r}$ belongs to a finite set. Since $\beta=\alpha_{k_0}$, $\beta$ belongs to a finite set, and we are done. 

\medskip

\noindent\textbf{Step 6}. In this step we prove (5). For any $1\leq j\leq r-1$, there are four cases.

\medskip

\noindent\textbf{Case 1}. $\alpha_j\not\in\Psi$. The equality follows from (2).

\medskip

\noindent\textbf{Case 2}. $\alpha_j=\beta$. The equality follows from (4.c). 

\medskip

\noindent\textbf{Case 3}. $\alpha_j=t\beta$ or $(t\beta)'$ for some $2\leq t\leq k-1$, then $\beta\in N^0\cap [0,\frac{1}{2}]^4$. By (4.a), $\beta=\alpha_{k_0}$ for some $1\leq k_0\leq r-1$. By (4.b), $2k_0+j=r$. By (4.c), $\beta(x_1)\geq\frac{1}{2}$. Thus $t=2$, $\beta(x_1)=\frac{1}{2}$, $\alpha_j=(2\beta)'$, and $\alpha_j(x_1)=0$. By (4.c) again,  $\beta(f)=2\beta(x_1)=1$, and
\begin{align*}
   \{\frac{ja_2}{r}\}+\{\frac{ja_3}{r}\}+\{\frac{ja_4}{r}\}&=3-2(\{\frac{k_0a_2}{r}\}+\{\frac{k_0a_3}{r}\}+\{\frac{k_0a_4}{r}\})\\
   &=3-2(\{\frac{k_0a_1}{r}\}+\frac{k_0}{r})=\frac{j}{r}+1.
\end{align*}
Since $\alpha_j(x_1)=0$ and $2a_1\equiv e\mod r$, $\{\frac{ja_1}{r}\}=0$, and $\{\frac{je}{r}\}=0$. Therefore, $$\sum_{i=1}^4\{\frac{ja_i}{r}\}=\{\frac{je}{r}\}+\frac{j}{r}+1,$$ and the equality holds.

\medskip

\noindent\textbf{Case 4}. $\alpha_j=\beta'\mod \Zz$. By (4.a), $\beta\equiv\alpha_{k_0}$ for some $1\leq k_0\leq r-1$. Thus $k_0+j=r$. If $\beta=\alpha_{k_0}$, then since $\alpha_j=\beta'$, $\alpha_{k_0}\in N^0\cap (0,1)^4$. By (4.c),
$$\sum_{i=1}^4\{\frac{ja_i}{r}\}=\sum_{i=1}^4\{\frac{(r-k_0)a_i}{r}\}=4-\sum_{i=1}^4\{\frac{k_0a_i}{r}\}=4-(\{\frac{k_0e}{r}\}+\frac{k_0}{r}+1)= \{\frac{je}{r}\}+\frac{j}{r}+1.$$ If $\beta\not=\alpha_{k_0}$, then $r\mid k_0 a_i$ for some $i\in \{1,2,3,4\}$. By (3), there are two sub-cases.

\medskip

\noindent\textbf{Case 4.1}. $\gcd(a_4,r)=\gcd(e,r)=2$ and $\gcd(a_1,r)=\gcd(a_2,r)=\gcd(a_3,r)=1$. Then $r\mid k_0 a_4$, and $k_0=\frac{r}{2}$. Hence $\beta$ belongs to a finite set, and we are done.

\medskip

\noindent\textbf{Case 4.2}. $\gcd(a_1,r)=\gcd(e,r)$ and $\gcd(a_2,r)=\gcd(a_3,r)=\gcd(a_4,r)=1$. In this case, $\{\frac{k_0a_1}{r}\}=0$, $\alpha_j(x_1)=\{\frac{ja_1}{r}\}=0$, and $\{\frac{je}{r}\}=0$. 

\medskip

\noindent\textbf{Case 4.2.1}. If $\alpha_{k_0}\in\Psi$, then since $\alpha_{k_0}\not=\beta$, $\alpha_{k_0}\equiv t\beta\mod\Zz^4$ for some $t\in\{k-1,\dots,2,-1,-2,\dots,-(k-1)\}$. Thus $r\mid (t-1)k_0$ for some $t\in\{k-1,\dots,2,-1,-2,\dots,-(k-1)\}$, hence $\frac{k_0}{r}$ belongs to a finite set. Since $\beta\equiv\alpha_{k_0}\mod\Zz^4$ and $\beta\in N^0$, $\beta$ belongs to a finite set, and we are done.

\medskip

\noindent\textbf{Case 4.2.2.2}. If $\alpha_{k_0}\not\in\Psi$, then by (2), 
$$\{\frac{k_0a_2}{r}\}+\{\frac{k_0a_3}{r}\}+\{\frac{k_0a_4}{r}\}=\frac{k_0}{r}+1.$$
Since $k_0+j=r$,
$$\{\frac{ja_1}{r}\}+\{\frac{ja_2}{r}\}+\{\frac{ja_3}{r}\}+\{\frac{ja_4}{r}\}=0+\frac{j}{r}+1=\{\frac{je}{r}\}+\frac{j}{r}+1,$$
and we are done.

\medskip

\noindent\textbf{Step 7}. In this step, we prove (6). By (3), (5), and the terminal lemma (Theorem \ref{thm: terminal lemma}), if $\gcd(a_1,r)=\gcd(e,r)\geq 2$, then $a_1\equiv e\mod r$. Since $2a_1\equiv e\mod r$, $a_1\equiv e\equiv 0\mod r$.
\end{proof}

\subsubsection{Odd type}

\begin{prop}\label{prop: odd case}
Notations and conditions as in Setting \ref{Setting: before terminal lem}. For each positive integer $k$, there exists a finite set $\Ii_k'$ depending only on $k$ satisfying the following. Suppose that $f$ is of odd type, and $\frac{1}{r}(a_1,a_2,a_3,a_4,e)\equiv\frac{1}{r}(1,\frac{r+2}{2},\frac{r-2}{2},2,2)\mod \Zz^5$ such that $4\mid r$. Then either $r\in\Ii_k'$ or $\bm{0}\not=\beta\in\Ii_k'$.
\end{prop}
\begin{proof}
We may assume that $r>4$, and $\beta\not\in\Ii_k$, where $\Ii_k$ is the set as in Theorem \ref{thm: non-cA case up to terminal lemma}. 

We consider 
$$\alpha_{r-2}=\frac{1}{r}(r-2,r-2,2,r-4).$$
There are two cases. 

\medskip

\noindent\textbf{Case 1}. $\alpha_{r-2}\not\in\Psi$. In this case, since $\alpha_{r-2}(x_2x_3x_4)=\frac{2r-4}{r}<\frac{2r-2}{r}=\alpha_{r-2}(x_1)+1$, by Theorem \ref{thm: non-cA case up to terminal lemma}(1), $\alpha_{r-2}(f)=2\alpha_{r-2}(x_1)-1=\frac{r-4}{r}$. Since $\alpha_{r-2}(x_1^2)=\alpha_{r-2}(x_2^2)=\frac{2r-4}{r}\not=\frac{r-4}{r}$, there exists a monomial $\bm{x}\in g\in\mm^3$ such that $\alpha_{r-2}(\bm{x})=\frac{r-4}{r}$. Since $\alpha_{r-2}(x_4)=\frac{r-4}{r}$, $\bm{x}=x_3^l$ for some integer $l\geq 3$. Thus $2l=r-4$, and $l=\frac{r-4}{2}$. Hence $$2\equiv e\equiv \alpha_1(x_3^l)=\frac{r-2}{2}\cdot\frac{r-4}{2}=2+\frac{r(r-6)}{4}\mod r.$$
Since $4\mid r$, $2+\frac{r(r-6)}{4}\equiv 2+\frac{r}{2}\not\equiv 2\mod r$, a contradiction.

\medskip

\noindent\textbf{Case 2}. $\alpha_{r-2}\in\Psi$. In this case, if $\alpha_{r-2}\in \Psi_2$, then $2\beta(x_1)\le 2(1-\alpha_{r-2}(x_1))=\frac{4}{r}<1$ which contradicts Theorem \ref{thm: non-cA case up to terminal lemma}(4.c). Thus $\alpha_{r-2}\in \Psi_1$. Since $\alpha_{r-2}(x_1)=\frac{r-2}{r}<1$, by Theorem \ref{thm: non-cA case up to terminal lemma}(4.c) again, $\alpha_{r-2}=\beta$. By Theorem \ref{thm: non-cA case up to terminal lemma}(4.b), $\frac{r-2}{r}<\frac{13}{14}$, so $r\leq 27$ and we are done.
\end{proof}

\subsubsection{cD-E type}

\begin{prop}\label{prop: cDE b case}
Notations and conditions as in Setting \ref{Setting: before terminal lem}. For each positive integer $k$, there exists a finite set $\Ii_k'$ depending only on $k$ satisfying the following. Suppose that $f$ is of cD-E type, and $\frac{1}{r}(a_1,a_2,a_3,a_4,e)\equiv\frac{1}{r}(a,-a,1,2a,2a)\mod \Zz^5$ for some integer $a$ such that $\gcd(a,r)=1$ and $2\mid r$. Then either $r\in\Ii_k'$ or $\bm{0}\not=\beta\in\Ii_k'$.
\end{prop}
\begin{proof}
We may assume that $r\geq 13$. We may also assume that $\beta\not\in\Ii_k$ where $\Ii_k$ is the set as in Theorem \ref{thm: non-cA case up to terminal lemma}.

Since $\gcd(a,r)=1$, there exists an integer $1\leq b\leq r-1$ such that $ba\equiv\frac{r+2}{2}<r$, and we have
$$\alpha_b=\frac{1}{r}(\frac{r+2}{2},\frac{r-2}{2},b,2).$$
There are three cases.

\medskip

\noindent\textbf{Case 1}. $\alpha_b\not\in\Psi$. In this case, since 
$$\alpha_b(x_2x_3x_4)=\frac{r+2}{2r}+\frac{b}{r}<\frac{r+2}{2r}+1=\alpha_b(x_1)+1,$$ by  Theorem \ref{thm: non-cA case up to terminal lemma}(1), $\alpha_b(f)=2\alpha_b(x_1)-1=\frac{2}{r}$. Thus there exists a monomial $\bm{x}\in\mm^3$ such that $\alpha_b(\bm{x})=\alpha_b(f)=\frac{2}{r}$, which is impossible.

\medskip

\noindent\textbf{Case 2}. $\alpha_b\in\Psi_2$. In this case, $2\beta(x_1)\le 2(1-\alpha_b(x_1))=\frac{r-2}{r}<1$ which contradicts Theorem \ref{thm: non-cA case up to terminal lemma}(4.c).

\medskip

\noindent\textbf{Case 3}. $\alpha_b\in\Psi_1$. In this case, since $\alpha(x_1)=\frac{r+2}{2r}<1$, by Theorem \ref{thm: non-cA case up to terminal lemma}(4.c), $\alpha_b=\beta$. In particular, $\Psi\cap N^0=\{\beta,\beta'\}=\{\frac{1}{r}(\frac{r+2}{2},\frac{r-2}{2},b,2),\frac{1}{r}(\frac{r-2}{2},\frac{r+2}{2},r-b,r-2)\}$

Since $\gcd(a,r)=1$, there exists an integer $1\leq c\leq r-1$ such that $ca\equiv\frac{r+4}{2}<r$, and we have
$$\alpha_c=\frac{1}{r}(\frac{r+4}{2},\frac{r-4}{2},c,4).$$
Since $\alpha_c \in N^0$, $\alpha_c(x_1)>\max\{\beta(x_1),\beta'(x_1)\}$, and $\Psi\cap N^0=\{\beta,\beta'\}$, we have $\alpha_c \not\in \Psi$. Since $$\alpha_c(x_2x_3x_4)=\frac{r+4}{2r}+\frac{c}{r}<\frac{r+4}{2r}+1=\alpha_c(x_1)+1,$$ by Theorem \ref{thm: non-cA case up to terminal lemma}(1), $\alpha_{c}(f)=2\alpha_c(x_1)-1=\frac{4}{r}$. Thus there exists a monomial $\bm{x}\in g\in\mm^3$ such that $\alpha_c(f)=\alpha_c(\bm{x})=\frac{4}{r}$. Since $r\geq 13$, $\alpha_c(x_2)=\frac{r-4}{2r}>\frac{4}{r}$ and $\alpha_c(x_4)=\frac{4}{r}$. Thus $\bm{x}=x_3^l$ for some integer $l\geq 3$, and $\alpha_c(x_3^l)=\frac{4}{r}$. Thus $lc=4$, hence $l=4$ and $c=1$. Thus $a\equiv\frac{r+4}{2}\mod r$ and $b\frac{r+4}{2}\equiv\frac{r+2}{2}\mod r$. Thus $r\mid 4b-2$. Since $2\mid r$ and $b<r$, either $r=4b-2$ or $3r=4b-2$. There are two sub-cases:

\medskip

\noindent\textbf{Case 3.1}. $r=4b-2$. In this case, we may assume that $b\geq 5$, otherwise $r\leq 18$ and we are done. Then $\alpha_b=\frac{1}{4b-2}(2b,2b-2,b,2)$ and $\alpha_c=\alpha_1=\frac{1}{4b-2}(2b+1,2b-3,1,4)$. We let $1\leq d\leq r-1$ be the unique positive integer such that $d\equiv 2b(b+1)\mod 4b-2$. Then $\alpha_d=\frac{1}{4b-2}(2b+2,2b-4,d,6)$. It is clear that $\alpha_d\not\in\{\beta,\beta'\}$. Since $2\alpha_d(x_1)=\frac{4b+4}{4b-2}>1$, by  Theorem \ref{thm: non-cA case up to terminal lemma}(1), $\alpha_d(f)=2\alpha_d(x_1)-1=\frac{6}{4b-2}$. Thus there exists a monomial $\bm{y}\in g\in\mm^3$ such that $\alpha_d(\bm{y})=\alpha_d(f)=\frac{6}{4b-2}$. Since $\alpha_d(x_2)=\frac{2b-4}{4b-2}\geq\frac{6}{4b-2}=\alpha_d(x_4)$, $\bm{y}=x_3^s$ for some positive integer $s$. Since
$$\alpha_d(x_3^s)=\alpha_d(f)\leq\alpha_d(x_3^l)$$
and
$$\alpha_1(x_3^l)=\alpha_1(f)\leq\alpha_d(x_3^s),$$
we have $s=l=4$. However,
$$\frac{ld}{4b-2}=\alpha_d(x_3^l)=\frac{6}{4b-2},$$
and $ld=6$, a contradiction.

\medskip

\noindent\textbf{Case 3.2}. $3r=4b-2$. In this case, we have $b=3s+2$ and $r=4s+2$ for some integer $s$. We may assume that $s\geq 4$, otherwise $r\leq 14$ and we are done. Then $\alpha_b=\beta=\frac{1}{4s+2}(2s+2,2s,3s+2,2)$ and $\alpha_c=\alpha_1=\frac{1}{4s+2}(2s+3,2s-1,1,4)$. We let $1\leq d\leq r-1$ be the unique positive integer such that $d\equiv -2s(s+2)\mod 4s+2$. Then $\alpha_d=\frac{1}{4s+2}(2s+4,2s-2,d,6)$. It is clear that $\alpha_d\not\in\{\beta,\beta'\}=\Psi\cap N^0$. Since $2\alpha_d(x_1)=\frac{4s+8}{4s+2}>1$, by Theorem \ref{thm: non-cA case up to terminal lemma}(1), $\alpha_d(f)=2\alpha_d(x_1)-1=\frac{6}{4s+2}$. Thus there exists a monomial $\bm{y}\in g\in\mm^3$ such that $\alpha_d(\bm{y})=\alpha_d(f)=\frac{6}{4s+2}$. Since $\alpha_d(x_2)=\frac{2s-2}{4s+2}\geq\frac{6}{4s+2}=\alpha_d(x_4)$, $\bm{y}=x_3^s$ for some positive integer $s$. Since
$$\alpha_d(x_3^s)=\alpha_d(f)\leq\alpha_d(x_3^l)$$
and
$$\alpha_1(x_3^l)=\alpha_1(f)\leq\alpha_d(x_3^s),$$
we have $s=l=4$. However,
$$\frac{ld}{4s+2}=\alpha_d(x_3^l)=\frac{6}{4s+2},$$
hence $ld=6$, a contradiction.
\end{proof}

\begin{prop}\label{prop: cDE c case}
Notations and conditions as in Setting \ref{Setting: before terminal lem}. For each positive integer $k$, there exists a finite set $\Ii_k'$ depending only on $k$ satisfying the following. Suppose that $f$ is of odd type, and $\frac{1}{r}(a_1,a_2,a_3,a_4,e)\equiv\frac{1}{r}(1,a,-a,2,2)\mod \Zz^5$ for some positive integer $a$ such that $\gcd(a,r)=1$ and $2\mid r$. Then either $r\in\Ii_k'$ or $\bm{0}\not=\beta\in\Ii_k'$.
\end{prop}
\begin{proof}
We may assume that $r\geq 13$ and $\beta\not\in\Ii_k$ where $\Ii_k$ is the set as in Theorem \ref{thm: non-cA case up to terminal lemma}, otherwise there is nothing left to prove.

Since $\gcd(a,r)=1$, we may let $1\leq b\leq r-1$ be the unique integer such that $b\equiv\frac{r+2}{2}a\mod r$. Then
$$\alpha_{\frac{r+2}{2}}=\frac{1}{r}(\frac{r+2}{2},b,r-b,2).$$
There are three cases.

\medskip

\noindent\textbf{Case 1}. $\alpha_{\frac{r+2}{2}}\not\in\Psi$. Since $2\alpha_{\frac{r+2}{2}}(x_1)=\frac{r+2}{r}>1$, by Theorem \ref{thm: non-cA case up to terminal lemma}(1), $\alpha_{\frac{r+2}{2}}(f)=2\alpha_{\frac{r+2}{2}}(x_1)-1=\frac{2}{r}$.  Thus there exists a monomial $\bm{x}\in g\in\mm^3$ such that $\alpha_{\frac{r+2}{2}}(\bm{x})=\alpha_{\frac{r+2}{2}}(f)=\frac{2}{r}$, which is impossible.

\medskip

\noindent\textbf{Case 2}.  $\alpha_{\frac{r+2}{2}}\in\Psi_2$. Then $2\beta(x_1)\le 2(1-\alpha_{\frac{r+2}{2}}(x_1))=\frac{r-2}{r}<1$, which contradicts Theorem \ref{thm: non-cA case up to terminal lemma}(4.c).

\medskip

\noindent\textbf{Case 3}. $\alpha_{\frac{r+2}{2}}\in\Psi_1$. Since $2\beta(x_1)\le 2(1-\alpha_{\frac{r+2}{2}}(x_1))=\frac{r-2}{r}<1$or $2\beta$. If $\alpha_{\frac{r+2}{2}}=2\beta$, then $2\beta(x_1)<1$, which contradicts Theorem \ref{thm: non-cA case up to terminal lemma}(4.c). Thus $\alpha_{\frac{r+2}{2}}=\beta$. In particular, $\Psi\cap N^0=\{\beta,\beta'\}=\{\frac{1}{r}(\frac{r+2}{2},b,r-b,2),\frac{1}{r}(\frac{r-2}{2},r-b,b,r-2)\}$.

We consider $\alpha_{\frac{r+4}{2}}=\frac{1}{r}(\frac{r+4}{2},c,r-c,4)$ where $c$ is the unique positive integer such that $c\equiv\frac{r+4}{2}a\mod r$. It is clear that $\alpha_{\frac{r+4}{2}}\not\in\{\beta,\beta'\}=\Psi$. Since $2\alpha_{\frac{r+4}{2}}(x_1)>1$, by Theorem \ref{thm: non-cA case up to terminal lemma}(1), $\alpha_{\frac{r+4}{2}}(f)=2\alpha_{\frac{r+4}{2}}(x_1)-1=\frac{4}{r}$. Thus there exists a monomial $\bm{x}\in g\in\mm^3$ such that $\alpha_{{\frac{r+4}{2}}}(\bm{x})=\frac{4}{r}$. Since $\alpha_{\frac{r+4}{2}}(x_4)=\frac{4}{r}$, $\alpha_{\frac{r+4}{2}}(x_2)=\frac{c}{r}$, and $\alpha_{\frac{r+4}{2}}(x_3)=\frac{r-c}{r}$, $c=1$ or $r-1$. Thus $1\equiv\frac{r+4}{2}a\mod r$ or $-1\equiv\frac{r+4}{2}a\mod r$, hence $\frac{r+2}{2}\equiv\frac{r+4}{2}b\mod r$ or $-\frac{r+2}{2}\equiv\frac{r+4}{2}b\mod r$. Thus $r\mid 4b+2$ or $r\mid 4b-2$. Since $b<r$ and $2\mid r$, there are four cases:

\medskip

\noindent\textbf{Case 3.1}. $r=4b-2$. We may assume that $b\geq 8$, otherwise $r\leq 30$ and we are done. In this case, $\beta=\alpha_{2b}=\frac{1}{4b-2}(2b,b,3b-2,2)$. Since $4b-2=r\mid 2(b-a)$ and $\gcd(a,r)=1$, we have $a=3b-1$, $\alpha_1=\frac{1}{4b-2}(1,3b-1,b-1,2)$, $c=1$, and $\alpha_{\frac{r+4}{2}}=\alpha_{2b+1}=\frac{1}{4b-2}(2b+1,1,4b-3,4)$.

We consider $\alpha_{2b+2}=\frac{1}{4b-2}(2b+2,3b,b-2,6)$. It is clear that $\alpha_{2b+2}\not\in\{\beta,\beta'\}=\Psi\cap N^0$. Then since $2\alpha_{2b+2}(x_1)=\frac{4b+4}{4b-2}>1$, by Theorem \ref{thm: non-cA case up to terminal lemma}(1), $\alpha_{2b+2}(f)=2\alpha_{2b+2}(x_1)-1=\frac{6}{4b-2}$. Thus there exits a monomial $\bm{y}\in g\in\mm^3$ such that $\alpha_{2b+2}(f)=\frac{6}{r}$. Since $b\geq 8$, $\alpha_{2b+2}(x_2)\geq\frac{6}{r}$, $\alpha_{2b+2}(x_3)\geq\frac{6}{r}$, and $\alpha_{2b+2}(x_4)=\frac{6}{r}$, a contradiction.

\medskip

\noindent\textbf{Case 3.2}. $3r=4b-2$. Then $b=3s+2$ and $r=4s+2$ for some positive integer $s$. We may assume that $s\geq 4$, otherwise $r\leq 14$ and we are done. Then $\alpha_{2s+2}=\beta=\frac{1}{4s+2}(2s+2,3s+2,s,2)$. Since $4s+2=r\mid 2(b-a)$ and $\gcd(a,r)=1$, we have $a=s+1$, $\alpha_1=\frac{1}{4s+2}(1,s+1,3s+1,2)$, $c=1$, and $\alpha_{\frac{r+4}{2}}=\alpha_{2s+3}=\frac{1}{4s+2}(2s+3,1,4s+1,4)$. We consider $\alpha_{2s+4}=\frac{1}{4s+2}(2s+4,s+2,3s,6)$. It is clear that $\alpha_{2s+4}\not\in\{\beta,\beta'\}=\Psi\cap N^0$. Then since $2\alpha_{2s+4}(x_1)=\frac{4s+8}{4s+2}>1$, by Theorem \ref{thm: non-cA case up to terminal lemma}(1), $\alpha_{2s+4}(f)=2\alpha_{2s+4}(x_1)-1=\frac{6}{4s+2}$. Thus there exits a monomial $\bm{y}\in g\in\mm^3$ such that $\alpha_{2s+4}(f)=\frac{6}{r}$. Since $s\geq 4$, $\alpha_{2s+4}(x_2)\geq\frac{6}{r}$, $\alpha_{2s+4}(x_3)\geq\frac{6}{r}$, and $\alpha_{2s+4}(x_4)=\frac{6}{r}$, a contradiction.

\medskip

\noindent\textbf{Case 3.3}. $r=4b+2$. We may assume that $b\geq 4$, otherwise $r\leq 14$ and we are done. In this case, $\beta=\alpha_{2b+2}=\frac{1}{4b+2}(2b+2,b,3b+2,2)$. Since $4b+2=r\mid 2(b-a)$ and $\gcd(a,r)=1$, we have $a=3b+1$, $\alpha_1=\frac{1}{4b+2}(1,3b+1,b+1,2)$, $c=r-1$, and $\alpha_{\frac{r+4}{2}}=\alpha_{2b+3}=\frac{1}{4b+2}(2b+3,4b+1,1,4)$.

We consider $\alpha_{2b+4}=\frac{1}{4b+2}(2b+4,3b,b+2,6)$. It is clear that $\alpha_{2b+4}\not\in\{\beta,\beta'\}=\Psi\cap N^0$. Then since $2\alpha_{2b+4}(x_1)=\frac{4b+8}{4b+2}>1$, by Theorem \ref{thm: non-cA case up to terminal lemma}(1), $\alpha_{2b+4}(f)=2\alpha_{2b+4}(x_1)-1=\frac{6}{4b+2}$. Thus there exits a monomial $\bm{y}\in g\in\mm^3$ such that $\alpha_{2b+4}(f)=\frac{6}{r}$. Since $b\geq 4$, $\alpha_{2b+4}(x_2)\geq\frac{6}{r}$, $\alpha_{2b+4}(x_3)\geq\frac{6}{r}$, and $\alpha_{2b+4}(x_4)=\frac{6}{r}$, a contradiction.

\medskip

\noindent\textbf{Case 3.2}. $3r=4b+2$. Then $b=3s+1$ and $r=4s+2$ for some positive integer $s$. We may assume that $s\geq 7$, otherwise $r\leq 26$ and we are done. Then $\alpha_{2s+2}=\beta=\frac{1}{4s+2}(2s+2,3s+1,s+1,2)$. Since $4s+2=r\mid 2(b-a)$ and $\gcd(a,r)=1$, we have $a=s$, $\alpha_1=\frac{1}{4s+2}(1,s,3s+2,2)$, $c=r-1$, and $\alpha_{\frac{r+4}{2}}=\alpha_{2s+3}=\frac{1}{4s+2}(2s+3,4s+1,1,4)$. We consider $\alpha_{2s+4}=\frac{1}{4s+2}(2s+4,s-1,3s+3,6)$. It is clear that $\alpha_{2s+4}\not\in\{\beta,\beta'\}=\Psi\cap N^0$. Then since $2\alpha_{2s+4}(x_1)=\frac{4s+8}{4s+2}>1$, by Theorem \ref{thm: non-cA case up to terminal lemma}(1), $\alpha_{2s+4}(f)=2\alpha_{2s+4}(x_1)-1=\frac{6}{4s+2}$. Thus there exits a monomial $\bm{y}\in g\in\mm^3$ such that $\alpha_{2s+4}(f)=\frac{6}{r}$. Since $s\geq 7$, $\alpha_{2s+4}(x_2)\geq\frac{6}{r}$, $\alpha_{2s+4}(x_3)\geq\frac{6}{r}$, and $\alpha_{2s+4}(x_4)=\frac{6}{r}$, a contradiction.
\end{proof}

\begin{prop}\label{prop: cDE d case}
Notations and conditions as in Setting \ref{Setting: before terminal lem}. For each positive integer $k$, there exists a finite set $\Ii_k'$ depending only on $k$ satisfying the following. Suppose that $f$ is of cD-E type, and $\frac{1}{r}(a_1,a_2,a_3,a_4,e)\equiv\frac{1}{r}(\frac{r-1}{2},\frac{r+1}{2},a,-a,-1)\mod \Zz^5$ for some integer positive $a$ such that $\gcd(a,r)=1$ and $2\nmid r$. Then either $r\in\Ii_k'$ or $\bm{0}\not=\beta\in\Ii_k'$.
\end{prop}
\begin{proof}
We may assume that $r\geq 13$ and $\beta\not\in\Ii_k$ where $\Ii_k$ is the set as in Theorem \ref{thm: non-cA case up to terminal lemma}, otherwise there is nothing left to prove. We may assume that $1\leq a\leq r-1$, then $\alpha_1=\frac{1}{r}(\frac{r-1}{2},\frac{r+1}{2},a,r-a)$, and $$\alpha_{r-1}=\frac{1}{r}(\frac{r+1}{2},\frac{r-1}{2},r-a,a).$$ There are three cases.

\medskip

\noindent\textbf{Case 1}. $\alpha_{r-1}\not\in\Psi$. Since $2\alpha_{r-1}(x_1)=\frac{r+1}{r}>1$, by Theorem \ref{thm: non-cA case up to terminal lemma}, $\alpha_{r-1}(f)=2\alpha_{r-1}(x_1)-1=\frac{1}{r}$. Thus there exists a monomial $\bm{x}\in g\in\mm^3$ such that $\alpha_{r-1}(\bm{x})=\frac{1}{r}$, which is impossible.

\medskip

\noindent\textbf{Case 2}. $\alpha_{r-1}\in\Psi_2$. Then $2\beta(x_1)\leq 2(1-\alpha_{r-1}(x_1))=\frac{r-1}{r}<1$, which contradicts Theorem \ref{thm: non-cA case up to terminal lemma}(4.c).

\medskip

\noindent\textbf{Case 3}. $\alpha_{r-1}\in\Psi_1$. Since $\gcd(\frac{r+1}{2},\frac{r-1}{2})=1$, $\alpha_{r-1}=\beta$. By Theorem \ref{thm: non-cA case up to terminal lemma}(4.b), $\frac{r-1}{r}<\frac{13}{14}$, hence $r\leq 13$, and we are done.
\end{proof}

\begin{prop}\label{prop: cDE e case}
Notations and conditions as in Setting \ref{Setting: before terminal lem}. For each positive integer $k$, there exists a finite set $\Ii_k'$ depending only on $k$ satisfying the following. Suppose that $f$ is of odd type, and $\frac{1}{r}(a_1,a_2,a_3,a_4,e)\equiv\frac{1}{r}(a,-a,2a,1,2a)\mod \Zz^5$ for some positive integer $a$ such that $\gcd(a,r)=1$ and $2\nmid r$. Then either $r\in\Ii_k'$ or $\bm{0}\not=\beta\in\Ii_k'$.
\end{prop}
\begin{proof}
We may assume that $r\geq 15$ and $\beta\not\in\Ii_k$ where $\Ii_k$ is the set as in Theorem \ref{thm: non-cA case up to terminal lemma}, otherwise there is nothing left to prove. 

Since $\gcd(a,r)=1$ and $2\nmid r$, there exists a unique positive integer $1\leq b\leq r-1$ such that $ba\equiv\frac{r+1}{2}\mod r$. Then
$$\alpha_{b}=\frac{1}{r}(\frac{r+1}{2},\frac{r-1}{2},1,b).$$
There are three cases.

\medskip

\noindent\textbf{Case 1}. $\alpha_b\not\in\Psi$.  Since $2\alpha_{b}(x_1)=\frac{r+1}{r}>1$, by Theorem \ref{thm: non-cA case up to terminal lemma}, $\alpha_{b}(f)=2\alpha_{b}(x_1)-1=\frac{1}{r}$. Thus there exists a monomial $\bm{x}\in g\in\mm^3$ such that $\alpha_{b}(\bm{x})=\frac{1}{r}$, which is impossible.

\medskip

\noindent\textbf{Case 2}. $\alpha_b\not\in\Psi_2$. Then $2\beta(x_1)\leq 2(1-\alpha_b)=\frac{r-1}{r}<1$, which contradicts Theorem \ref{thm: non-cA case up to terminal lemma}(4.c).

\medskip

\noindent\textbf{Case 3}. $\alpha_{b}\in\Psi_1$. Since $\gcd(\frac{r+1}{2},\frac{r-1}{2})=1$, $\alpha_{b}=\beta$. Thus $\Psi\cap N^0=\{\beta,\beta'\}=\{\frac{1}{r}(\frac{r+1}{2},\frac{r-1}{2},1,b),\frac{1}{r}(\frac{r-1}{2},\frac{r+1}{2},r-1,r-b)\}$.

We let $c$ and $d$ be the unique integers such that $ca\equiv\frac{r+3}{2} \mod r$ and $da\equiv\frac{r+5}{2} \mod r$. Then
$$\alpha_c=\frac{1}{r}(\frac{r+3}{2},\frac{r-3}{2},3,c)$$
and
$$\alpha_d=\frac{1}{r}(\frac{r+5}{2},\frac{r-5}{2},5,d).$$
Since $r\geq 13$, it is clear that $\alpha_c\not\in\Psi$ and $\alpha_{d}\not\in\Psi$. Since $2\alpha_c(x_1)>1$ and $2\alpha_d(x_1)>1$, by Theorem \ref{thm: non-cA case up to terminal lemma}(1), $\alpha_c(f)=2\alpha_c(x_1)-1=\frac{3}{r}$ and $\alpha_d(f)=2\alpha_{d}(x_1)-1=\frac{5}{r}$. Thus there exist monomials $\bm{x},\bm{y}\in g\in\mm^3$ such that $\alpha_c(\bm{x})=\alpha_c(f)=\frac{3}{r}$ and $\alpha_d(\bm{y})=\alpha_d(f)=\frac{5}{r}$. Since $r\geq 15$, $\alpha_c(x_2)\geq\frac{6}{r}$, $\alpha_d(x_2)\geq\frac{5}{r}$, $\alpha_c(x_3)=\frac{3}{r}$ and $\alpha_d(x_3)=\frac{5}{r}$. Thus $\bm{x}=x_4^l$ and $\bm{y}=x_4^s$ for some $l,s\geq 3$. Moreover, since
$$\alpha_c(x_4^l)=\alpha_c(f)=\frac{3}{r}\leq\alpha_c(x_4^s)$$
and
$$\alpha_d(x_4^s)=\alpha_d(f)=\frac{5}{r}\le \alpha_d(x_4^l),$$
we have $l=s$. Since $\alpha_c(x_4)=\frac{c}{r}$, and $\alpha_d(x_4)=\frac{d}{r}$, we have $lc=3$ and $ld=5$, which contradicts $l\ge 3$.
\end{proof}

\begin{prop}\label{prop: cDE f case}
Notations and conditions as in Setting \ref{Setting: before terminal lem}. For each positive integer $k$, there exists a finite set $\Ii_k'$ depending only on $k$ satisfying the following. Suppose that $f$ is of cD-E type, and $\frac{1}{r}(a_1,a_2,a_3,a_4,e)\equiv\frac{1}{r}(1,a,-a,2,2)\mod \Zz^5$ for some positive integer $a$ such that $\gcd(a,r)=1$ and $2\nmid r$. Then either $r\in\Ii_k'$ or $\bm{0}\not=\beta\in\Ii_k'$.
\end{prop}
\begin{proof}
We may assume that $r\geq 15$ and $\beta\not\in\Ii_k$ where $\Ii_k$ is the set as in Theorem \ref{thm: non-cA case up to terminal lemma}, otherwise there is nothing left to prove. We may assume that $1\leq a\leq r-1$.

Since $\gcd(a,r)=1$ and $2\nmid r$, there exists a unique positive integer $1\leq b\leq r-1$ such that $b\equiv\frac{r+1}{2}a\mod r$. Then
$$\alpha_{\frac{r+1}{2}}=\frac{1}{r}(\frac{r+1}{2},b,r-b,1).$$
There are three cases.

\medskip

\noindent\textbf{Case 1}. $\alpha_{\frac{r+1}{2}}\not\in\Psi$. In this case, since $2\alpha_{\frac{r+1}{2}}(x_1)=\frac{r+1}{r}>1$, by Theorem \ref{thm: non-cA case up to terminal lemma}(1), $\alpha_{\frac{r+1}{2}}(f)=2\alpha_{\frac{r+1}{2}}(x_1)-1=\frac{1}{r}$. Thus there exists a monomial $\bm{x}\in g\in\mm^3$ such that $\alpha_{\frac{r+1}{2}}(\bm{x})=\frac{1}{r}$, which is impossible.

\medskip

\noindent\textbf{Case 2}. $\alpha_{\frac{r+1}{2}}\in\Psi_2$. Then $2\beta(x_1)\le 2(1-\alpha_{\frac{r+1}{2}})=\frac{r-1}{r}<1$, which contradicts Theorem \ref{thm: non-cA case up to terminal lemma}(4.c).

\medskip

\noindent\textbf{Case 3}. $\alpha_{\frac{r+1}{2}}\in\Psi_1$. Since $\alpha_{\frac{r+1}{2}}(x_4)=\frac{1}{r}$, $\alpha_{\frac{r+1}{2}}=\beta$. In particular, $\Psi\cap N^0=\{\beta,\beta'\}=\{\frac{1}{r}(\frac{r+1}{2},b,r-b,1),\frac{1}{r}(\frac{r-1}{2},r-b,b,r-1)\}$.

Let $c$ be the unique integer such that $c\equiv\frac{r+3}{2}a\mod r$. Then
$$\alpha_{\frac{r+3}{2}}=\frac{1}{r}(\frac{r+3}{2},c,r-c,3).$$
Since $r\geq 13$, it is clear that $\alpha_c\not\in\Psi$. Since $2\alpha_c(x_1)>1$, by Theorem \ref{thm: non-cA case up to terminal lemma}(1), $\alpha_c(f)=2\alpha_c(x_1)-1=\frac{3}{r}$. Thus there exists a monomial $\bm{x}\in g\in\mm^3$ such that $\alpha_{\frac{r+3}{2}}(\bm{x})=\alpha_c(f)=\frac{3}{r}$. Since $\alpha_{\frac{r+3}{2}}(x_2)=\frac{c}{r},\alpha_{\frac{r+3}{2}}(x_3)=\frac{r-c}{r}$, and $\alpha_{\frac{r+3}{2}}(x_4)=\frac{3}{r}$, $c=1$ or $r-1$. Thus $1\equiv\frac{r+3}{2}a\mod r$ or $-1\equiv\frac{r+3}{2}a\mod r$, so $\frac{r+1}{2}\equiv\frac{r+3}{2}b \mod r$ or $-\frac{r+1}{2}\equiv\frac{r+3}{2}b \mod r$. Thus $r\mid 3b-1$ or $r\mid 3b+1$. Since $b<r$, there are four cases.

\medskip

\noindent\textbf{Case 3.1}. $r=3b-1$. Since $2\nmid r$, $b=2s$ and $r=6s-1$ for some positive integer $s$. We may assume $s\geq 4$, otherwise $r\leq 17$ and we are done. Then $\alpha_{\frac{r+1}{2}}=\alpha_{3s}=\frac{1}{6s-1}(3s,2s,4s-1,1)=\beta$. Since $6s-1=r\mid 2b-a$ and $\gcd(a,r)=1$, we have $a=4s$, $\alpha_1=\frac{1}{6s-1}(1,4s,2s-1,2)$, and $\alpha_{\frac{r+3}{2}}=\alpha_{3s+1}=\frac{1}{6s-1}(3s+1,1,6s-2,3)$. We consider
$$\alpha_{3s+2}=\frac{1}{6s-1}(3s+2,4s+1,2s-2,5),$$
then it is clear that $\alpha_{3s+2}\not\in\Psi$. Since $2\alpha_{3s+2}(x_1)>1$, by Theorem \ref{thm: non-cA case up to terminal lemma}(1), $\alpha_{3s+2}(f)=2\alpha_{3s+2}(x_1)-1=\frac{5}{r}$. Thus there exists a monomial $\bm{x}\in g\in\mm^3$ such that $\alpha_{3s+2}(f)=\alpha(\bm{x})=\frac{5}{r}$. Since $s\geq 4$, $\alpha_{3s+2}(x_2)\geq\frac{5}{r},$ $\alpha_{3s+2}(x_3)\geq\frac{5}{r},$ and $\alpha_{3s+2}(x_4)=\frac{5}{r}$, a contradiction.

\medskip

\noindent\textbf{Case 3.2}. $2r=3b-1$. Since $2\nmid r$, $b=4s+1$ and $r=6s+1$ for some positive integer $s$. We may assume $s\geq 2$, otherwise $r\leq 7$ and we are done. Then $\alpha_{\frac{r+1}{2}}=\alpha_{3s+1}=\frac{1}{6s+1}(3s+1,4s+1,2s,1)=\beta$. Since $6s+1=r\mid 2b-a$ and $\gcd(a,r)=1$, we have $a=2s+1$, $\alpha_1=\frac{1}{6s+1}(1,2s+1,4s,2)$, and $\alpha_{\frac{r+3}{2}}=\alpha_{3s+2}=\frac{1}{6s+1}(3s+2,1,6s,3)$. We consider
$$\alpha_{3s+3}=\frac{1}{6s+1}(3s+3,2s+2,4s-1,5),$$
then it is clear that $\alpha_{3s+3}\not\in\Psi$. Since $2\alpha_{3s+3}(x_1)>1$, by Theorem \ref{thm: non-cA case up to terminal lemma}(1), $\alpha_{3s+3}(f)=2\alpha_{3s+3}(x_1)-1=\frac{5}{r}$. Thus there exists a monomial $\bm{x}\in g\in\mm^3$ such that $\alpha_{3s+3}(f)=\alpha(\bm{x})=\frac{5}{r}$. Since $s\geq 2$, $\alpha_{3s+3}(x_2)\geq\frac{5}{r},$ $\alpha_{3s+3}(x_3)\geq\frac{5}{r},$ and $\alpha_{3s+3}(x_4)=\frac{5}{r}$, a contradiction.

\medskip

\noindent\textbf{Case 3.3}. $r=3b+1$. Since $2\nmid r$, $b=2s$ and $r=6s+1$ for some positive integer $s$. We may assume $s\geq 2$, otherwise $r\leq 7$ and we are done. Then $\alpha_{\frac{r+1}{2}}=\alpha_{3s+1}=\frac{1}{6s+1}(3s+1,2s,4s+1,1)=\beta$. Since $6s+1=r\mid 2b-a$ and $\gcd(a,r)=1$, we have $a=4s$, $\alpha_1=\frac{1}{6s+1}(1,4s,2s+1,2)$, and $\alpha_{\frac{r+3}{2}}=\alpha_{3s+2}=\frac{1}{6s+1}(3s+2,6s,1,3)$. We consider
$$\alpha_{3s+3}=\frac{1}{6s+1}(3s+3,4s-1,2s+2,5),$$
then it is clear that $\alpha_{3s+3}\not\in\Psi$. Since $2\alpha_{3s+3}(x_1)>1$, by Theorem \ref{thm: non-cA case up to terminal lemma}(1), $\alpha_{3s+3}(f)=2\alpha_{3s+3}(x_1)-1=\frac{5}{r}$. Thus there exists a monomial $\bm{x}\in g\in\mm^3$ such that $\alpha_{3s+3}(f)=\alpha(\bm{x})=\frac{5}{r}$. Since $s\geq 2$, $\alpha_{3s+3}(x_2)\geq\frac{5}{r},$ $\alpha_{3s+3}(x_3)\geq\frac{5}{r},$ and $\alpha_{3s+3}(x_4)=\frac{5}{r}$, a contradiction.

\medskip

\noindent\textbf{Case 3.4}. $2r=3b+1$. Since $2\nmid r$, $b=4s-1$ and $r=6s-1$ for some positive integer $s$. We may assume $s\geq 4$, otherwise $r\leq 17$ and we are done. Then $\alpha_{\frac{r+1}{2}}=\alpha_{3s}=\frac{1}{6s-1}(3s,4s-1,2s,1)=\beta$. Since $6s-1=r\mid 2b-a$ and $\gcd(a,r)=1$, we have $a=2s-1$, $\alpha_1=\frac{1}{6s-1}(1,2s-1,4s,2)$, and $\alpha_{\frac{r+3}{2}}=\alpha_{3s+1}=\frac{1}{6s-1}(3s+1,6s-2,1,3)$. We consider
$$\alpha_{3s+2}=\frac{1}{6s-1}(3s+2,2s-2,4s+1,5),$$
then it is clear that $\alpha_{3s+2}\not\in\Psi$. Since $2\alpha_{3s+2}(x_1)>1$, by Theorem \ref{thm: non-cA case up to terminal lemma}(1), $\alpha_{3s+2}(f)=2\alpha_{3s+3}(x_1)-1=\frac{5}{r}$. Thus there exists a monomial $\bm{x}\in g\in\mm^3$ such that $\alpha_{3s+2}(f)=\alpha(\bm{x})=\frac{5}{r}$. Since $s\geq 4$, $\alpha_{3s+2}(x_2)\geq\frac{5}{r},$ $\alpha_{3s+2}(x_3)\geq\frac{5}{r},$ and $\alpha_{3s+2}(x_4)=\frac{5}{r}$, a contradiction.
\end{proof}

\begin{lem}\label{lem: cDE a case lemma}
Notations and conditions as in Setting \ref{Setting: before terminal lem}. For each positive integer $k$, there exists a finite set $\Ii_k'$ depending only on $k$ satisfying the following. Suppose that $f$ is of cD-E type, and $\frac{1}{r}(a_1,a_2,a_3,a_4,e)\equiv\frac{1}{r}(0,a,-a,1,0)\mod \Zz^5$ for some positive integer $a$ such that $\gcd(a,r)=1$. 

Then for any $1\leq j\leq r-1$ such that $\alpha_j\not\equiv t\beta\mod\Zz^4$ for any $1\leq t\leq k-1$, $\alpha_j(g)=1$.
\end{lem}
\begin{proof}
Since $\alpha_j(g)\equiv \alpha_j(x_1)=0\mod \Zz$, $\alpha_j(g)\in\Zz_{\ge 0}$ for any $j$. Since $\alpha_j(x_i)\not=0$ for any $1\leq j\leq r-1$ and $i\in\{2,3,4\}$, $\alpha_j(g)\in\Zz_{\ge 1}$.

For any $1\leq j\leq r-1$ such that $\alpha_j\not\equiv t\beta\mod\Zz^4$ for any $1\leq t\leq k-1$, we let $\gamma_j:=\alpha_j+(\lceil\frac{\alpha_j(g)}{2}\rceil,0,0,0)$. Then $\gamma_j\not\in\{\beta,2\beta,\dots,(k-1)\beta\}$. By Setting \ref{Setting: before terminal lem}(3.b.ii), $$\lceil\frac{\alpha_j(g)}{2}\rceil+\frac{r+j}{r}=\gamma_j(x_1x_2x_3x_4)>\gamma_j(f)+1=\min\{\alpha_j(g),2\lceil\frac{\alpha_j(g)}{2}\rceil\}+1=\alpha_j(g)+1.$$
Thus $\alpha_j(g)\leq\lceil\frac{\alpha_j(g)}{2}\rceil$, hence $\alpha_j(g)=1$.
\end{proof}

\begin{prop}\label{prop: cDE a case}
Notations and conditions as in Setting \ref{Setting: before terminal lem}. For each positive integer $k$, there exists a finite set $\Ii_k'$ depending only on $k$ satisfying the following. Suppose that $f$ is of odd type, and $\frac{1}{r}(a_1,a_2,a_3,a_4,e)\equiv\frac{1}{r}(0,a,-a,1,0)\mod \Zz^5$ for some positive integer $a$ such that $\gcd(a,r)=1$. Then either $r\in\Ii_k'$ or $\bm{0}\not=\beta\in\Ii_k'$.
\end{prop}
\begin{proof}
We may assume that $r\geq 6k+2$ otherwise there is nothing left to prove. We define $b_1,\dots,b_k$ in the following way: for any $1\leq j\leq k$, let $1\leq b_j\leq r-1$ be the unique integer such that $ab_j\equiv \lceil\frac{r}{2}\rceil-j \mod r$. Then $$\alpha_{b_j}=\frac{1}{r}(0,\lceil\frac{r}{2}\rceil-j,\lfloor\frac{r}{2}\rfloor+j,b_j).$$
and
$$\alpha_{r-b_j}=\frac{1}{r}(0,\lfloor\frac{r}{2}\rfloor+j,\lceil\frac{r}{2}\rceil-j,r-b_j).$$
for any $1\leq j\leq k$. By the pigeonhole principle, there exists $j_0\in\{1,2,\dots,k\}$ such that $\alpha_{b_{j_0}}\not\equiv t\beta$ and $\alpha_{r-b_{j_0}}\not\equiv t\beta$ for any $1\leq t\leq k-1$. Since $r\geq 6k+2$, $\lceil\frac{r}{2}\rceil-j_0>\frac{r}{3}$. Thus $\alpha_{b_{j_0}}(x_2)>\frac{1}{3}$, $\alpha_{b_{j_0}}(x_3)>\frac{1}{3}$, $\alpha_{r-b_{j_0}}(x_2)>\frac{1}{3}$, and $\alpha_{r-b_{j_0}}(x_3)>\frac{1}{3}$. By Lemma \ref{lem: cDE a case lemma}, $\alpha_{b_{j_0}}(g)=1$ and $\alpha_{r-b_{j_0}}(g)=1$. Thus $\frac{b_{j_0}}{r}=\alpha_{b_{j_0}}(x_4)\leq \frac{1}{3}$ and $\frac{r-b_{j_0}}{r}=\alpha_{r-b_{j_0}}(x_4)\leq \frac{1}{3}$, which is impossible.
\end{proof}

\begin{proof}[Proof of Theorem \ref{thm: non cA type beta finite}] 

We may assume $r\ge 3$. By Theorem \ref{thm: non-cA case up to terminal lemma}(3,5), $a_i,e$ satisfy the conditions of terminal lemma. 

If $f$ is of odd type, then $a_1\not\equiv a_2 \mod r$ and $2a_1\equiv 2a_2\equiv e \mod r$. Thus by the terminal lemma (Theorem \ref{thm: terminal lemma}) and Theorem \ref{thm: non-cA case up to terminal lemma}(6), possibly interchanging $a_1,a_2$, $\frac{1}{r}(a_1,a_2,a_3,a_4,e)\equiv \frac{1}{r}(1,\frac{r+2}{2},\frac{r-2}{2},2,2) \mod \Zz^5$ and $4\mid r$. Now the theorem follows from Proposition \ref{prop: odd case}.

If $f$ is of cD-E type, then $2a_1\equiv e \mod r$. Thus by the terminal lemma (Theorem \ref{thm: terminal lemma}), possibly interchanging $a_1,a_2$ or $a_3,a_4$, one of the following holds.
\begin{enumerate}
    \item $\frac{1}{r}(a_1,a_2,a_3,a_4,e)\equiv\frac{1}{r}(a,-a,1,2a,2a) \mod \Zz^5$ for some positive integer $a$ such that $\gcd(a,r)=1$ and $2\mid r$.
    \item $\frac{1}{r}(a_1,a_2,a_3,a_4,e)=\frac{1}{r}(1,a,-a,2,2) \mod \Zz^5$ for some positive integer $a$ such that $\gcd(a,r)=1$ and $2\mid r$.
    
    \item $\frac{1}{r}(a_1,a_2,a_3,a_4,e)\equiv\frac{1}{r}(\frac{r-1}{2},\frac{r+1}{2},a,-a,-1) \mod \Zz^5$ for some positive integer $a$ such that $\gcd(a,r)=1$ and $2\nmid r$.
    
    \item $\frac{1}{r}(a_1,a_2,a_3,a_4,e)\equiv\frac{1}{r}(a,-a,2a,1,2a) \mod \Zz^5$ for some positive integer $a$ such that $\gcd(a,r)=1$ and $2\nmid r$.
    
    \item $\frac{1}{r}(a_1,a_2,a_3,a_4,e)\equiv\frac{1}{r}(1,a,-a,2,2) \mod \Zz^5$ for some positive integer $a$ such that $\gcd(a,r)=1$ and $2\nmid r$.
    
    \item $\frac{1}{r}(a_1,a_2,a_3,a_4,e)\equiv\frac{1}{r}(0,a,-a,1,0) \mod \Zz^5$ for some positive integer $a$ such that $\gcd(a,r)=1$.
\end{enumerate}
Now the theorem follows from Propositions \ref{prop: cDE b case}, \ref{prop: cDE c case}, \ref{prop: cDE d case}, \ref{prop: cDE f case}, and \ref{prop: cDE a case}. 
\end{proof}

\begin{proof}[Proof of Theorem \ref{thm: note beta finite}]
It follows from Theorems \ref{thm: cA type beta finite} and \ref{thm: non cA type beta finite}.
\end{proof}

\end{document}